%%%%%%%%%%%%%%%%%%%%%%%%%%%%%%%%%%%%%%%%%%%%%%%%%% 
%	 last version of April 8, 
%%%%%%%%%%%%%%%%%%%%%%%%%%%%%%%
% begin.tex

%=========================================================================%
% load begin.tex only once, but keep count to match \bye commands
%=========================================================================%

\ifx\begin\undefined\else\global\advance\srcdepth by
1\expandafter \fi

\def\begin{}
\newcount\srcdepth
\srcdepth=1

\outer\def\bye{\global\advance\srcdepth by -1
  \ifnum\srcdepth=0
    \def\endcmd{\vfill\eject\nopagenumbers\par\vfill\supereject\end}
  \else\def\endcmd{}\fi
  \endcmd
}

%=========================================================================%
% initialize TeX
%=========================================================================%

\magnification=\magstephalf
%\magnification=\magstep 1
\baselineskip=13pt
%baselineskip=12pt
\hsize = 5.5truein
\hoffset = 0.5truein
\vsize = 8.5truein
\voffset = 0.2truein
\emergencystretch = 0.05\hsize

\overfullrule=0pt

\newif\ifblackboardbold

% comment out the following line if AMS msbm fonts aren't available
\blackboardboldtrue

%=========================================================================%
% select fonts
%=========================================================================%

\font\sectionfont=cmbx12

% Establish AMS blackboard bold fonts without using amssym.def, amssym.tex

\newfam\bboldfam
\ifblackboardbold
\font\tenbbold=msbm10
\font\sevenbbold=msbm7
\font\fivebbold=msbm5
\textfont\bboldfam=\tenbbold
\scriptfont\bboldfam=\sevenbbold
\scriptscriptfont\bboldfam=\fivebbold
\def\bbold{\fam\bboldfam\tenbbold}
\else
\def\bbold{\bf}
\fi

%=========================================================================%
% font size-changing command ("A Beginner's Book of TeX" p35, p275)
%=========================================================================%

\font\Arm=cmr8
\font\Ai=cmmi8
\font\Asy=cmsy8
\font\Abf=cmbx8
\font\Brm=cmr6
\font\Bi=cmmi6
\font\Bsy=cmsy6
\font\Bbf=cmbx6
\font\Crm=cmr5
\font\Ci=cmmi5
\font\Csy=cmsy5
\font\Cbf=cmbx5

\ifblackboardbold
\font\Abbold=msbm10 at 8pt
\font\Bbbold=msbm7 at 6pt
\font\Cbbold=msbm5
\fi

\def\smallmath{%
\textfont0=\Arm \scriptfont0=\Brm \scriptscriptfont0=\Crm
\textfont1=\Ai \scriptfont1=\Bi \scriptscriptfont1=\Ci
\textfont2=\Asy \scriptfont2=\Bsy \scriptscriptfont2=\Csy
\textfont\bffam=\Abf \scriptfont\bffam=\Bbf \scriptscriptfont\bffam=\Cbf
\def\rm{\fam0\Arm}\def\mit{\fam1}\def\oldstyle{\fam1\Ai}%
\def\bf{\fam\bffam\Abf}%
\ifblackboardbold
\textfont\bboldfam=\Abbold
\scriptfont\bboldfam=\Bbbold
\scriptscriptfont\bboldfam=\Cbbold
\def\bbold{\fam\bboldfam\Abbold}%
\fi
}

%=========================================================================%
% single-pass symbolic theorem labeling
%=========================================================================%

% Because this is a single-pass mechanism with no .aux file, forward
% references need to be declared in advance:

%   \forward{thm:main}{Theorem}{1.1}

% This is also the mechanism for "timely" declaration of labels, which
% will usually be buried within the corresponding theorem macros.
% A warning is issued if a label redeclaration is inconsistent, allowing
% forward references to be manually fixed.

%   \ref{thm:main} produces "Theorem~1.1"
%   \refs{thm:main} produces "Theorems~1.1"
%   \refn{thm:main} produces "1.1"

% Some TeX adapted from "The Advanced TeXbook" by David Salomon, chapter 9.

% Implementers: The code for \forward is subtle. Its second argument must
% be provided literally, e.g. "Theorem" rather that "\capitalize{theorem}".
% Its third argument must either be literal or a macro that expands
% directly to a literal, e.g. "\edef\numtoks{\number\proccount}".
% This use of \edef cannot be replaced by \def, which defers expansion.
% Failure to follow these rules will cause spurious warnings that forward
% references are inconsistent, when they are in fact consistent after
% expansion. Note the "Towers of Palo Alto" recreational math problem
% involving the iterated use of \expandafter to expand the first argument
% to \forwardsub before calling it.

\newlinechar=`@
\def\forwardmsg#1#2#3{\immediate\write16{@*!*!*!* forward reference should
be: @\noexpand\forward{#1}{#2}{#3}@}}
\def\nodefmsg#1{\immediate\write16{@*!*!*!* #1 is an undefined reference@}}

\def\forwardsub#1#2{\def\newref{{#2}{#1}}}

\def\forward#1#2#3{%
\expandafter\expandafter\expandafter\forwardsub\expandafter{#3}{#2}
\expandafter\ifx\csname#1\endcsname\relax\else%
\expandafter\ifx\csname#1\endcsname\newref\else%
\forwardmsg{#1}{#2}{#3}\fi\fi%
\expandafter\let\csname#1\endcsname\newref}

\def\firstarg#1{\expandafter\argone #1}\def\argone#1#2{#1}
\def\secondarg#1{\expandafter\argtwo #1}\def\argtwo#1#2{#2}

\def\ref#1{\expandafter\ifx\csname#1\endcsname\relax
  {\nodefmsg{#1}\bf`#1'}\else
  \expandafter\firstarg\csname#1\endcsname
  ~\expandafter\secondarg\csname#1\endcsname\fi}

\def\refs#1{\expandafter\ifx\csname#1\endcsname\relax
  {\nodefmsg{#1}\bf`#1'}\else
  \expandafter\firstarg\csname #1\endcsname
  s~\expandafter\secondarg\csname#1\endcsname\fi}

\def\refn#1{\expandafter\ifx\csname#1\endcsname\relax
  {\nodefmsg{#1}\bf`#1'}\else
  \expandafter\secondarg\csname #1\endcsname\fi}

%=========================================================================%
% widow control
%=========================================================================%

% usage:
% \widow{.2} % start new page if <.2 page left

\def\widow#1{\vskip 0pt plus#1\vsize\goodbreak\vskip 0pt plus-#1\vsize}

%=========================================================================%
% sections and theorems
%=========================================================================%

% use \showlabels or \showlabelsabove to display section and theorem labels

\def\marginlabel#1{}

\def\showlabelsabove{
\font\labelfont=cmss10 at 6pt
\def\marginlabel##1{\rlap{\smash{\raise 10pt\hbox{\labelfont##1}}}}
}

\newcount\seccount
\newcount\proccount
\seccount=0
\proccount=0

\def\stdskip{\vskip 9pt plus3pt minus 3pt}
\def\stdbreak{\par\removelastskip\penalty-100\stdskip}

\def\proof{\stdbreak\noindent{\sl Proof. }}

\def\qed{\vrule height 1.2ex width .9ex depth .1ex}

\def\Box{
  \ifmmode\eqno\qed
  \else\ifvmode\removelastskip\line{\hfil\qed}
  \else\unskip\quad\hskip-\hsize
    \hbox{}\hskip\hsize minus 1em\qed\par
  \fi\stdbreak\fi}

\def\references{
  \removelastskip
  \widow{.05}
  \vskip 24pt plus 6pt minus 6 pt
  \leftline{\sectionfont References}
  \nobreak\stdskip\noindent}

\def\ifempty#1#2\endB{\ifx#1\endA}
\def\makeref#1#2#3{\ifempty#1\endA\endB\else\forward{#1}{#2}{#3}\fi}

\outer\def\section#1 #2\par{
  \removelastskip
  \global\advance\seccount by 1
  \global\proccount=0\relax
                \edef\numtoks{\number\seccount}
  \makeref{#1}{Section}{\numtoks}
  \widow{.05}
  \vskip 24pt plus 6pt minus 6 pt
  \message{#2}
  \leftline{\marginlabel{#1}\sectionfont\numtoks\quad #2}
  \nobreak\stdskip}

\def\proclamation#1#2{
  \outer\expandafter\def\csname#1\endcsname##1 ##2\par{
  \stdbreak
  \advance\proccount by 1
  \edef\numtoks{\number\seccount.\number\proccount}
  \makeref{##1}{#2}{\numtoks}
  \noindent{\marginlabel{##1}\bf #2 \numtoks\enspace}
  {\sl##2\par}
  \stdbreak}}

\def\othernumbered#1#2{
  \outer\expandafter\def\csname#1\endcsname##1{
  \stdbreak
  \advance\proccount by 1
  \edef\numtoks{\number\seccount.\number\proccount}
  \makeref{##1}{#2}{\numtoks}
  \noindent{\marginlabel{##1}\bf #2 \numtoks\enspace}}}

\proclamation{definition}{Definition}
\proclamation{lemma}{Lemma}
\proclamation{proposition}{Proposition}
\proclamation{theorem}{Theorem}
\proclamation{corollary}{Corollary}
\proclamation{conjecture}{Conjecture}

\othernumbered{example}{Example}
\othernumbered{remark}{Remark}
\othernumbered{construction}{Construction}

%=========================================================================%
% enable postscript illustrations using epsf.tex
%=========================================================================%

% Usage:
% \draw{70}{fig}{} % draw fig.eps at 70% scale
% \draw{999}{fig}{} % draw fig.eps scaled to width of page

% Optional third argument can be multiple calls to \figtext; see below.
% More generally, the third argument is read in vertical mode, with the
% reference point at the lower left corner of the eps picture, whose
% dimensions are contained in the dimen registers \drawx and \drawy.
% This enables using TeX to generate the text that goes with the picture.
% To request that the picture be widened to respect the added text, 
% examine and modify the dimen registers \ngap, \egap, \sgap, \wgap.
% This is done automatically by the \figtext macro.

% These macros rely on "epsf.tex" which is the lowest level interface
% available for including encapsulated Postscript files in TeX documents.
% Rather that manually reading the .eps file to compute the nominal size,
% the \epsfbox macro is called twice, and two of its internal registers
% are examined after the first call. A major change to epsf.tex (unlikely)
% will require changes here. 

\input epsf

\newcount\figcount
\figcount=0
\newbox\drawing
\newcount\drawbp
\newdimen\drawx
\newdimen\drawy
\newdimen\ngap
\newdimen\sgap
\newdimen\wgap
\newdimen\egap

\def\drawbox#1#2#3{\vbox{
  \setbox\drawing=\vbox{\offinterlineskip\epsfbox{#2.eps}\kern 0pt}
  \drawbp=\epsfurx
  \advance\drawbp by-\epsfllx\relax
  \multiply\drawbp by #1
  \divide\drawbp by 100
  \drawx=\drawbp truebp
  \ifdim\drawx>\hsize\drawx=\hsize\fi
  \epsfxsize=\drawx
  \setbox\drawing=\vbox{\offinterlineskip\epsfbox{#2.eps}\kern 0pt}
  \drawx=\wd\drawing
  \drawy=\ht\drawing
  \ngap=0pt \sgap=0pt \wgap=0pt \egap=0pt 
  \setbox0=\vbox{\offinterlineskip
    \box\drawing \ifgridlines\drawgrid\drawx\drawy\fi #3}
  \kern\ngap\hbox{\kern\wgap\box0\kern\egap}\kern\sgap}}

\def\draw#1#2#3{
  \setbox\drawing=\drawbox{#1}{#2}{#3}
  \advance\figcount by 1
  \goodbreak
  \midinsert
  \centerline{\ifgridlines\boxgrid\drawing\fi\box\drawing}
  \smallskip
  \vbox{\offinterlineskip
    \centerline{Figure~\number\figcount}
    \smash{\marginlabel{#2}}}
  \endinsert}

\def\nextfigtoks{%
  \advance\figcount by 1%
  \edef\numtoks{\number\figcount}%
  \advance\figcount by -1}

\newif\ifgridlines
\newbox\figtbox
\newbox\figgbox
\newdimen\figtx
\newdimen\figty

\newdimen\bwd
\bwd=2sp % 2sp (1/32768") is smallest visible width for Textures

\def\hline#1{\vbox{\smash{\hbox to #1{\leaders\hrule height \bwd\hfil}}}}

\def\vline#1{\hbox to 0pt{%
  \hss\vbox to #1{\leaders\vrule width \bwd\vfil}\hss}}

\def\clap#1{\hbox to 0pt{\hss#1\hss}}
\def\vclap#1{\vbox to 0pt{\offinterlineskip\vss#1\vss}}

\def\hstutter#1#2{\hbox{%
  \setbox0=\hbox{#1}%
  \hbox to #2\wd0{\leaders\box0\hfil}}}

\def\vstutter#1#2{\vbox{
  \setbox0=\vbox{\offinterlineskip #1}
  \dp0=0pt
  \vbox to #2\ht0{\leaders\box0\vfil}}}

\def\crosshairs#1#2{
  \dimen1=.002\drawx
  \dimen2=.002\drawy
  \ifdim\dimen1<\dimen2\dimen3\dimen1\else\dimen3\dimen2\fi
  \setbox1=\vclap{\vline{2\dimen3}}
  \setbox2=\clap{\hline{2\dimen3}}
  \setbox3=\hstutter{\kern\dimen1\box1}{4}
  \setbox4=\vstutter{\kern\dimen2\box2}{4}
  \setbox1=\vclap{\vline{4\dimen3}}
  \setbox2=\clap{\hline{4\dimen3}}
  \setbox5=\clap{\copy1\hstutter{\box3\kern\dimen1\box1}{6}}
  \setbox6=\vclap{\copy2\vstutter{\box4\kern\dimen2\box2}{6}}
  \setbox1=\vbox{\offinterlineskip\box5\box6}
  \smash{\vbox to #2{\hbox to #1{\hss\box1}\vss}}}

\def\boxgrid#1{\rlap{\vbox{\offinterlineskip
  \setbox0=\hline{\wd#1}
  \setbox1=\vline{\ht#1}
  \smash{\vbox to \ht#1{\offinterlineskip\copy0\vfil\box0}}
  \smash{\vbox{\hbox to \wd#1{\copy1\hfil\box1}}}}}}

\def\drawgrid#1#2{\vbox{\offinterlineskip
  \dimen0=\drawx
  \dimen1=\drawy
  \divide\dimen0 by 10
  \divide\dimen1 by 10
  \setbox0=\hline\drawx
  \setbox1=\vline\drawy
  \smash{\vbox{\offinterlineskip
    \copy0\vstutter{\kern\dimen1\box0}{10}}}
  \smash{\hbox{\copy1\hstutter{\kern\dimen0\box1}{10}}}}}

\def\figtext#1#2#3#4#5{
  \setbox\figtbox=\hbox{#5}
  \dp\figtbox=0pt
  \figtx=-#3\wd\figtbox \figty=-#4\ht\figtbox
  \advance\figtx by #1\drawx \advance\figty by #2\drawy
  \dimen0=\figtx \advance\dimen0 by\wd\figtbox \advance\dimen0 by-\drawx
  \ifdim\dimen0>\egap\global\egap=\dimen0\fi
  \dimen0=\figty \advance\dimen0 by\ht\figtbox \advance\dimen0 by-\drawy
  \ifdim\dimen0>\ngap\global\ngap=\dimen0\fi
  \dimen0=-\figtx
  \ifdim\dimen0>\wgap\global\wgap=\dimen0\fi
  \dimen0=-\figty
  \ifdim\dimen0>\sgap\global\sgap=\dimen0\fi
  \smash{\rlap{\vbox{\offinterlineskip
    \hbox{\hbox to \figtx{}\ifgridlines\boxgrid\figtbox\fi\box\figtbox}
    \vbox to \figty{}
    \ifgridlines\crosshairs{#1\drawx}{#2\drawy}\fi
    \kern 0pt}}}}

% macros to add space to text on specified sides

\def\hpad#1#2#3{\hbox{\kern #1\hbox{#3}\kern #2}}
\def\vpad#1#2#3{\setbox0=\hbox{#3}\dp0=0pt\vbox{\kern #1\box0\kern #2}}

% macro to give one text string the apparent height of another

% macro to center one text string over another

\def\stack#1#2#3{\vbox{\offinterlineskip
  \setbox2=\hbox{#2}
  \setbox3=\hbox{#3}
  \dimen0=\ifdim\wd2>\wd3\wd2\else\wd3\fi
  \hbox to \dimen0{\hss\box2\hss}
  \kern #1
  \hbox to \dimen0{\hss\box3\hss}}}

% macros to hide size of trailing exponents

\def\hexp#1{%
  \setbox0=\hbox{${}^{#1}$}%
  \hbox to .5\wd0{\box0\hss}}

%=========================================================================%
% macros for matrices and arrows
%=========================================================================%

% typical usage:
%   \rightarrowmat{2pt}{4pt}{d & bd \cr \!-c & 0 \cr 0 & -ac \cr}

\def\bmatrix#1#2{{\smallmath\left[\vcenter{\halign
  {&\kern#1\hfil$##\mathstrut$\kern#1\cr#2}}\right]}}

\def\rightarrowmat#1#2#3{
  \setbox1=\hbox{\kern#2$\bmatrix{#1}{#3}$\kern#2}
  \,\vbox{\offinterlineskip\hbox to\wd1{\hfil\copy1\hfil}
    \kern 3pt\hbox to\wd1{\rightarrowfill}}\,}

\def\leftarrowmat#1#2#3{
  \setbox1=\hbox{\kern#2$\bmatrix{#1}{#3}$\kern#2}
  \,\vbox{\offinterlineskip\hbox to\wd1{\hfil\copy1\hfil}
    \kern 3pt\hbox to\wd1{\leftarrowfill}}\,}

\def\rightarrowbox#1#2{
  \setbox1=\hbox{\kern#1\hbox{\smallmath #2}\kern#1}
  \,\vbox{\offinterlineskip\hbox to\wd1{\hfil\copy1\hfil}
    \kern 3pt\hbox to\wd1{\rightarrowfill}}\,}

\def\leftarrowbox#1#2{
  \setbox1=\hbox{\kern#1\hbox{\smallmath #2}\kern#1}
  \,\vbox{\offinterlineskip\hbox to\wd1{\hfil\copy1\hfil}
    \kern 3pt\hbox to\wd1{\leftarrowfill}}\,}

%=========================================================================%
% quire macros for preview mode and making booklets
%=========================================================================%

% \legalbooklet{20} makes a booklet from legal paper in landscape
% orientation, where "20" is the page count. To preview, give a negative
% pagecount. Either print using the legal duplex option on a modern laser
% printer, or struggle to simulate this effect manually. Bind using a long
% reach stapler.

% \preview squeezes two pages side by side in landscape orientation. It
% is not suitable for printing, but ideal for previewing on a two page
% monitor.

% \twoup squeezes two pages onto letter paper in landscape mode,
% suitable for printing.

% Each of these macros calls the file "quire.tex"

\def\bookletdims{
  \hsize=5.25truein
  \vsize=7truein
}

\def\legalbooklet#1{
  \input quire
  \bookletdims
  \htotal=7.0truein
  \vtotal=8.5truein
  % below computed from above
  \hoffset=\htotal
  \advance\hoffset by -\hsize
  \divide\hoffset by 2
  \voffset=\vtotal
  \advance\voffset by -\vsize
  \divide\voffset by 2
  \advance\voffset by -.0625truein
  \shhtotal=2\htotal
  % below doesn't need to change
  \horigin=0.0truein
  \vorigin=0.0truein
  \shstaplewidth=0.01pt
  \shstaplelength=0.66truein
  \shthickness=0pt
  \shoutline=0pt
  \shcrop=0pt
  \shvoffset=-1.0truein
  \ifnum#1>0\quire{#1}\else\qtwopages\fi
}

\def\preview{
  \input quire
  \bookletdims
  \hoffset=0.1truein
  \vtotal=8.5truein
  \shhtotal=14truein
  % below computed from above
  \voffset=\vtotal
  \advance\voffset by -\vsize
  \divide\voffset by 2
  \advance\voffset by -.0625truein
  \htotal=2\hoffset
  \advance\htotal by \hsize
  % below doesn't need to change
  \horigin=0.0truein
  \vorigin=0.0truein
  \shstaplewidth=0.5pt
  \shstaplelength=0.5\vtotal
  \shthickness=0pt
  \shoutline=0pt
  \shcrop=0pt
  \shvoffset=-1.0truein
  \qtwopages
}

\def\twoup{
  \input quire
  \hsize=4.79452truein % 5.25/1.095
  \vsize=7truein
  \vtotal=8.5truein
  \shhtotal=11truein
  % below computed from above
  \hoffset=-2\hsize
  \advance\hoffset by \shhtotal
  \divide\hoffset by 6
  \voffset=\vtotal
  \advance\voffset by -\vsize
  \divide\voffset by 2
  \advance\voffset by -12truept
  \htotal=2\hoffset
  \advance\htotal by \hsize
  % below doesn't need to change
  \horigin=0.0truein
  \vorigin=0.0truein
  \shstaplewidth=0.01pt
  \shstaplelength=0pt
  \shthickness=0pt
  \shoutline=0pt
  \shcrop=0pt
  \shvoffset=-1.0truein
  \qtwopages
}

%=========================================================================%
% timestamp (adapted from eplain.tex)
%=========================================================================%

\newcount\countA
\newcount\countB
\newcount\countC

\def\monthname{\begingroup
  \ifcase\number\month
    \or January\or February\or March\or April\or May\or June\or
    July\or August\or September\or October\or November\or December\fi
\endgroup}

\def\dayname{\begingroup
  \countA=\number\day
  \countB=\number\year
  \advance\countA by 0 % adjust after each leap day
  \advance\countA by \ifcase\month\or
    0\or 31\or 59\or 90\or 120\or 151\or
    181\or 212\or 243\or 273\or 304\or 334\fi
  \advance\countB by -1995
  \multiply\countB by 365
  \advance\countA by \countB
  \countB=\countA
  \divide\countB by 7
  \multiply\countB by 7
  \advance\countA by -\countB
  \advance\countA by 1
  \ifcase\countA\or Sunday\or Monday\or Tuesday\or Wednesday\or
    Thursday\or Friday\or Saturday\fi
\endgroup}

\def\timename{\begingroup
   \countA = \time
   \divide\countA by 60
   \countB = \countA
   \countC = \time
   \multiply\countA by 60
   \advance\countC by -\countA
   \ifnum\countC<10\toks1={0}\else\toks1={}\fi
   \ifnum\countB<12 \toks0={\sevenrm AM}
     \else\toks0={\sevenrm PM}\advance\countB by -12\fi
   \relax\ifnum\countB=0\countB=12\fi
   \hbox{\the\countB:\the\toks1 \the\countC \thinspace \the\toks0}
\endgroup}

\def\timestamp{\dayname, \the\day\ \monthname\ \the\year, \timename}

%==========================================================================
% macros (specific to this paper)
%==========================================================================

% surround with $ $ if not already in math mode
\def\enma#1{{\ifmmode#1\else$#1$\fi}}

%\showlabels
\showlabelsabove
\input diagrams.tex
%\vsize=12truecm
%\hsize=22truecm
\overfullrule=0pt
% Gothic fonts from AMSTeX 
\font\tengoth=eufm10  \font\fivegoth=eufm5
\font\sevengoth=eufm7
\newfam\gothfam  \scriptscriptfont\gothfam=\fivegoth 
\textfont\gothfam=\tengoth \scriptfont\gothfam=\sevengoth
\def\goth{\fam\gothfam\tengoth}
%
% Bold italic fonts 
\font\tenbi=cmmib10  \font\fivebi=cmmib5
\font\sevenbi=cmmib7
\newfam\bifam  \scriptscriptfont\bifam=\fivebi 
\textfont\bifam=\tenbi \scriptfont\bifam=\sevenbi

\font\hd=cmbx10 scaled\magstep1
\def \fix#1 {{\hfill\break \bf (( #1 ))\hfill\break}}
\def \Box {\hfill\hbox{}\nobreak \vrule width 1.6mm height 1.6mm
depth 0mm  \par \goodbreak \smallskip}
\def \reg {\mathop{\rm regularity}}
\def \coker {\mathop{\rm coker}}
\def \ker {\mathop{\rm ker}}
\def \im {\mathop{\rm im}}

\def \dim{{\rm dim}}

\def \iso {\cong}
\def \tensor {\otimes}

\def \Hom {{\mathop{\rm Hom}\nolimits}}
\def \hom {{\mathop{\rm Hom}\nolimits}}
\def \Ext {{\rm Ext}}
\def \ext{{\rm Ext}}
\def \Tor {{\rm Tor}}
\def \tor{{\rm Tor}}
\def \Sym {{\mathop{\rm Sym}\nolimits}}
\def \sym{{\mathop{\rm Sym}\nolimits}}

\def \lin{{\rm lin}}

\def \th {{^{\rm th}}}

\def \AA {{\bf A}}

\def \B {{\cal B}}
\def \C {{\cal C}}

\def \F {{\cal F}}
\def \FF {{\bf F}}

\def \GG {{\bf G}}

\def \H {{\rm H}}

\def \L {{\cal L}}
\def \LL {{\bf L}}
\def \MM{{\bf M}}

\def \O {{\cal O}}
\def \P {{\bf P}}
\def \PP {{\bf P}}
\def \RR {{\bf R}}
\def \TT {{\bf T}}
\def \Z {{\bf Z}}

\def \gm {{\goth m}}

%%%%%%%%%%%%
\forward {notation}{Section}{1}
\forward {intro BGG}{Section}{2}
\forward {powers example}{Section}{5}
%\forward {reg and ex}{Section}{4}
\forward {linear part section}{Section}{3}
%\forward {resolutions}{Section}{6}
\forward {tate}{Section}{4}
\forward {beilinson}{Section}{6}
%\forward {socle res}{Section}{9}
%
%
\forward {basic correspondence}{Proposition}{2.1}
\forward {exactness criterion}{Corollary}{2.4}
\forward {linear exactness}{Corollary}{2.5}

\forward {eg1}{Example}{3.3}
\forward {linear part and tor}{Theorem}{3.4}
\forward {degenerate double complex}{Lemma}{3.5}

\forward {linear part 1}{Corollary}{3.6}
\forward {linear dominance}{Theorem}{3.1}
%\forward {regularity}{Corollary}{3.7}
%\forward {deformation 1}{Lemma}{3.8}

\forward {ecoh-thm1}{Theorem}{2.6}
\forward {BGG-result}{Corollary}{2.7}

\forward {reciprocity}{Theorem}{3.7}
%\forward {singular vectors of a standard module}{Corollary}{6.5}

\forward {sheaf cohomology}{Theorem}{4.1}
\forward {comp of shf coho}{Corollary}{4.2}
%\forward {ecoh-thm}{Theorem}{7.4}
\forward {lin cplx}{Corollary}{4.4}
\forward {elliptic quartic}{Example}{4.6}
\forward {hom of omega}{Proposition}{5.6}
\forward {beilinson-thm}{Theorem}{6.1}
\forward {beilinson cor 1}{Corollary}{6.3}
\forward {beilinson cor 2}{Corollary}{6.2}
\forward {beilinson-thm 2}{Theorem}{6.4}
\forward {row bound}{Lemma}{7.3}
%\forward {socle absence}{Theorem}{9.2}

\rightline {April 13, 2000}
\bigskip

\centerline{\hd Sheaf Cohomology} 
\centerline{\hd and}
\centerline{\hd Free Resolutions over Exterior Algebras}
\medskip 
\centerline {\bf David Eisenbud and Frank-Olaf Schreyer
\footnote{$^{*}$}{\rm Both authors are grateful to the NSF for 
partial support during the preparation of this paper. 
The second author wishes to thank MSRI for its hospitality.}
\footnote{}{\rm AMS Classification. Primary: 14F05, 14Q20, 16E05}
}
\bigskip
\bigskip

{\narrower
\noindent{\bf Abstract:} 
In this paper we derive an explicit version of the
Bernstein-Gel'fand-Gel'fand (BGG) correspondence between bounded
complexes of coherent sheaves on projective space and minimal doubly
infinite free resolutions over its ``Koszul dual'' exterior algebra.
This leads to an efficient method for machine computation of the
cohomology of sheaves. Among the facts about the BGG correspondence
that we derive is that taking homology of a complex of sheaves
corresponds to taking the ``linear part'' of a resolution over the
exterior algebra.

Using these results we give a constructive proof of the existence of a
Beilinson monad for a sheaf on projective space. The explicitness of
our version allows us to to prove two conjectures about
the morphisms in the monad.

Along the way we prove a number of results about minimal free resolutions
over an exterior algebra. For example, we show that such resolutions are
eventually dominated by their ``linear parts" in the sense that
erasing all terms of  degree $>1$ in the complex yields a new complex
which is eventually exact.

\bigskip}

Let V be a finite dimensional vector space over a field $K$, and let
$W=V^*$ be the dual space. In this paper we will study complexes and
resolutions over the exterior algebra $E=\wedge V$ and their relation
to modules over $S=\sym W$ and sheaves on projective space $\P(W)$.

Perhaps the most important known results in this area are:
\smallskip

\item{$\bullet$} The Bernstein-Gel'fand-Gel'fand  (BGG) 
correspondence [1978] (see also
Gel'fand [1984],
Buchweitz [1987], 
and Gel'fand-Manin [1996]), 
usually
stated as an equivalence between the derived category of 
bounded complexes of coherent sheaves on $\P(W)$ and
the derived category of bounded complexes of finitely
generated graded modules over $E$.

\item{$\bullet$} Beilinson's theorem [1978],
which gives, for each sheaf $\F$ on projective space, a complex
$$
\dots\rTo 
\oplus_{j=0}^n\H^j(\F(e-j))\otimes_K\Omega_{\P^n}^{j-e}(j-e)
\rTo\dots
$$
called the {\it Beilinson Monad\/} whose homology is precisely
$\F$ and whose terms depend only on the cohomology of a few twists of
$\F$. 
\smallskip

\noindent In this paper we re-examine the BGG correspondence. Its 
essential content is a functor $\RR$ from complexes of graded $S$-modules
to complexes of graded $E$-modules, and its adjoint $\LL$.
For example, if $M=\oplus_i M_i$ is a graded $S$-module (regarded
as a complex with just one term) then as a
bigraded
$E$-module $\RR(M)=\Hom_K(E,M)$, with differential 
$\Hom_K(E,M_i)\to \Hom_K(E,M_{i+1})$ defined from the multiplication map
on $M$. Similarly, for a graded $E$-module $P$, we have
$\LL(P)=S\otimes_KP$. In fact (\ref{basic correspondence}) $\RR$ is an
equivalence from the category of graded
$S$-modules and 
to the category of linear complexes of free $E$-modules; here {\it
linear\/} means essentially that the maps are represented by matrices
of linear forms. A similar statement holds for $\LL$.

We show that finitely generated modules $M$ go to left-bounded complexes
that are eventually exact on the right, and characterize the point at
which  exactness begins as the Castelnuovo-Mumford regularity of $M$. 
A strong form of this is \ref{reciprocity}, of which the following is
a part:

\proclaim Reciprocity Theorem.
If $M$ is a graded $S$-module and $P$ is a
graded
$E$-module, then $\RR(M)$ is an injective resolution of $P$ if and only if
$\LL(P)$ is a free resolution of $M$.

Let  $\F$ be a coherent sheaf on projective space and with $r>>0$ take
$M=\oplus_{d} \H^0(\F(d))$. The results above show that the complex
$
\RR(M_{\geq r})
$
associated to the $r^\th$ truncation of $M$
is acyclic for $r>>0$. If we take a minimal free resolution of the kernel
of the first term in this complex, we obtain a doubly infinite
exact free complex, independent of $r$, which we call
the {\it Tate resolution\/}  $\TT(\F)$:
$$
\TT(\F): \cdots \to T^{r-1}\to T^r=\Hom_K(E,M_r)\to \Hom_K(E,
M_{r+1})\to\cdots
$$
It was
first studied in Gel'fand [1984]. We show (\ref{sheaf cohomology}) that
the $e^\th$ term of the Tate resolution is
$T^e(\F)=\oplus_j\Hom_K(E, \H^j(\F(e-j))$; that is it is made
from the cohomology of the twists of $\F$.
This leads to a new algorithm for computing sheaf
cohomology. We have  programmed
this method in the computer algebra system Macaulay2 of Grayson
and Stillman  [{\tt http://www.math.uiuc.edu/Macaulay2/}]. In
some cases it gives the fastest known computation of the
cohomology. 

As a further application
of \ref{sheaf cohomology}, we show that if
$P=\ker\bigl[T^0(\F)\to T^1(\F)\bigr]$ then $\LL(P)$ is
a {\it linear monad} for $\F$: it is a
linear free complex of finitely generated $S$-modules
such that, for every $i$,
$$
\H^i(\LL(P))=\oplus_{j\geq -i}\H^i(\F(-j));
$$
and in particular the only homology of the sheafification of
$\LL(P)$ is $\F$ itself  (\ref{lin
cplx}).

One of the main results of this paper is a new proof of 
Beilinson's Theorem (\ref{beilinson-thm}), which clarifies its connection
of the BGG-correspondence. 
Beilinson's
original paper [1978] sketches a proof that leads easily to a weak
form of the result, the ``Beilinson spectral sequence'', which
determines the sheaf $\F$ only up to filtration. This version is
explained in the book of Okonek, Schneider, and Spindler
[1980]. Kapranov [1988] and Ancona and Ottaviani [1989] have given a full proof working in the derived category along the
lines indicated by Beilinson. However the use of the derived category 
makes it difficult to compute the 
Beilinson's monad for explicit sheaves effectively.

Our proof of Beilinson's Theorem
leads to a practical construction of the monad and to  new results about
its structure. Because of the explicitness of our construction, we can
answer some questions about the maps in the Beilinson monad: There are
natural  candidates for the linear components of the maps in the monad
for a sheaf $\F$; and given such 
a monad, there are natural candidates
for most of the maps in the  monad of $\F(1)$. Our techniques
allow us to prove that these natural candidates 
really do occur (\ref{beilinson cor 2} and \ref{beilinson cor 1}).
Beilinson [1978] also proved the existence of
a different monad for a sheaf $\F$, 
using the sheaves $\O_\P(i)$ for $0\leq i\leq v$
in place of the $\Omega^i(i)$. The second monad can be deduced easily from the the first monad (\ref{beilinson-thm 2}).

A remarkable feature of the theory of resolutions
over the exterior algebra, not visible for
the corresponding theory over a polynomial ring, is that the linear
terms of any  resolution eventually predominate.
To state this precisely, we introduce the
{\it linear part\/} of a free complex $\FF$ over $S$ or $E$. The linear
part 
is the complex obtained from $\FF$ by taking a minimal free
complex $\GG$ homotopic to $\FF$, and then erasing all terms of
absolute degree $>1$ from the matrices representing
the differentials of $\GG$. In fact taking the linear
part is functorial in a suitable sense: 
under the BGG correspondence it corresponds to the  homology
functor (\ref{linear part and tor}). 
Just as the homology of a complex is simpler than the complex,
one can often compute the linear part of a complex even when the
complex itself is mysterious. 

Of course free resolutions may have maps with no linear
terms at all, that is, with linear part equal
to zero. And they can have infinitely many maps with nonlinear
terms unavoidably present (this is even the case for periodic
resolutions; see Avramov and Eisenbud [2000]). But the linear
terms eventually predominate in the following
sense:

\proclaim Theorem \refn{linear dominance}. If\/ $\FF$ is the free
resolution of  a finitely generated  module 
over the exterior algebra $E$
then the linear part of\/ $\FF$ is eventually exact.

This predominance can take arbitrarily long to assert itself: the
resolution of the millionth syzygy of the residue field of
$E$ has a million linear maps follows by a map with
linear part 0, and linear dominance happens only at the 
million and first term. In the case of a resolution
of a monomial ideal, however, Herzog and R\"omer [1999] have
shown that the linear part becomes exact after at most
$v$ steps. It would be interesting to know more results of
this sort.

In the last section we illustrate the material of the paper with some
examples. 

Much of the elementary material of this paper could be done for an
arbitrary pair of homogeneous Koszul algebras (in the sense of Priddy
[1970]) in place of the pair of algebras $S, E$. We use a tiny bit of
this for the pair $(E,S)$.  See Buchweitz [1987] for a sketch of the
general case and a statement of general conditions under which the BGG
correspondence holds.  Versions of Beilinson's theorem have been
established for some other varieties through work of Swan [1985],
Kapranov [1988,1989], and Orlov [1992].  Yet other derived category
equivalences have been pursued under the rubric of ``tilting" (see
Happel [1988]).

This paper owes much to the experiments we were able to make
using the computer algebra system Macaulay2 of Grayson and Stillman,
and we would like to thank them for their help and patience with
this project. We are also grateful to Luchezar Avramov for
getting us interested in resolutions over exterior algebras.

\section{notation} Notation and Background

Throughout this paper we write $K$ for a fixed field, and $V, W$ for
dual vector spaces of finite dimension $v$ over $K$.  We give the
elements of $W$ degree 1, so that the elements of $V$ have degree
$-1$.  We write $E=\wedge V$ and $S=\Sym(W)$ for the exterior and
symmetric algebras; these algebras are graded by their {\it internal
degrees\/} whereby $\Sym_i(W)$ has degree $i$ and $\wedge^jV$ has
degree $-j$.  We think of $E$ as $Ext_S^\bullet(K,K)$ and $S$ as
$Ext_E^\bullet(K,K)$. 

We will always write the index indicating the degree of
a homogeneous component of
a graded module as subscripts. Thus if $M=\oplus M_i$ is a graded module over
$E$ or $S$, then $M_i$ denotes the component of degree $i$.
We let $M(a)$ be the shifted module, so that
$M(a)_b=M_{a+b}$. 
We write complexes cohomologically, with
upper indices and differentials of degree $+1$. Thus if 
$$
\FF:\qquad \dots \rTo F^i\rTo F^{i+1} \dots ,
$$ 
is a complex, 
then $F^i$ denotes the term of cohomological degree $i$.
We write $\FF\langle a\rangle$ for the complex whose
term of cohomological degree $j$ is $F^{a+j}$. 

We will write 
$\omega_S=S\otimes_K\wedge^vW$ 
for the module associated to the canonical bundle of $\P(W)$;
note that $\wedge^vW$ is a vector space concentrated in 
degree $v$, so that $\omega_S$ is noncanonically isomorphic
to $S(-v)$. Similarly, we set 
$\omega_E:=\Hom_K(E,K)=E\otimes_K\wedge^vW$, which is 
noncanonically isomorphic to $E(-v)$. 
It is easy to check that for any graded vector space $D$ we have 
$\Hom_K(E,D)\cong \omega_E\otimes_K D$
as left $E$-modules.
For any $E$-module $P$, we set $P^*:=\Hom_K(P,K)$.

We often use the fact that the exterior algebra is Gorenstein and 
finite dimensional over $K$,
which follows from the fact that $\Hom_K(E,K)\iso E$ as above.
As a consequence, the dual of any exact sequence is exact and the notions
free module, injective module, and projective module coincide.

We also use the notion of Castelnuovo-Mumford regularity.
The most convenient
definition for our purposes is that the Castelnuovo-Mumford
regularity of a graded $S$-module $M=\oplus_iM_i$ is the smallest integer
$r$ such that the truncation $M_{\geq r}=\oplus_{i\geq r}M_i$
is generated by $M_r$ and has
a {\it linear free resolution}---that is, all the maps in
its free resolution are represented by matrices of linear forms.
See for example
Eisenbud-Goto [1984] or Eisenbud [1995] 
for a discussion. The regularity of a sheaf $\F$ on projective space
(equal to the regularity of $\oplus_d \H^0(\F(d))$ if this module
is finitely generated) can also be expressed as the minimal
$r$ for which $\H^i(\F(r-i))=0$ for all $i>0$.

A free complex over $E$ or a graded free complex over $S$ is called {\it
minimal\/} if all its maps can be represented by matrices with entries in
the appropriate maximal ideal. For example, any linear complex is
minimal.

\section{intro BGG} The Bernstein-Gel'fand-Gel'fand Correspondence

In this section we will give a brief exposition of the main
idea of Bernstein-Gel'fand-Gel'fand [1978]: a construction of a 
pair of adjoint functors between the categories of complexes over $E$
and over $S$. However, we avoid the peculiar convention, used in the
original, according to which the differentials of complexes over $E$
were not homomorphisms of $E$-modules.

Let $e_i$ and $x_i$ be dual bases of $V$ and $W$ respectively, so that
$\sum_i x_i\otimes e_i\in W\otimes_K V$ corresponds to the identity
element under the isomorphism $W\otimes_K V = \Hom_K(W,W)$.
Let $A$ and $B$ be vector spaces. Giving a map 
$A\otimes_K W\rTo^\alpha B$ 
is the same as giving a map $A\rTo^{\alpha'}B\otimes_K V$ (where the
tensor products are taken over $K$). For example, given $\alpha$ we
set $\alpha'(a) = \sum_i e_i\otimes \alpha(x_i\otimes a)$.  

We begin with a special case that will play a central role. We 
regard a graded $S$-module $M=\oplus M_d$ 
as a complex with only one term, in cohomological degree 0, and 
define $\RR(M)$ to be the complex
$$
\eqalign{
\dots \rTo^\phi
\Hom_K(E, M_d)&\rTo^\phi
\Hom_K(E, M_{d+1})\rTo^\phi
 \dots  \cr
\phi:\ \alpha&\mapsto \bigl[ e\mapsto \sum_i
x_i\alpha(e_ie)\bigr]. }
$$
Here the term $\Hom_K(E,M_d)$ has cohomological index $d$, and a map
$\phi\in\Hom_K(E, M_{d})$ has degree $t$ if it factors through the 
projection from $E$ onto $E_{d-t}$. Note that the complex
$\RR(M)$ is {\it linear\/} in a strong sense:
the $d^\th$ free module $\Hom_K(E, M_d)\iso \omega_E\otimes M_d$ 
has socle in degree $d$;
in particular all the maps are represented by matrices of linear
forms. 

\proposition{basic correspondence} The functor $\RR$
is an equivalence between the category of graded left $S$-modules
and the category of linear free complexes over $E$ (those
for which the $d^\th$ free module has socle in degree $d$.)

\proof A collection of maps $\mu_d: W\otimes_K M_d\rTo M_{d+1}$ 
defines a module structure on the graded 
vector space $\oplus M_d$ if and only if it satisfies
a commutativity and associativity condition
expressed by saying that,
for each $d$, the composition of the multiplication maps
$$
W\otimes_K (W\otimes_K M_d)\rTo W\otimes_K M_{d+1}
\rTo M_{d+2}
$$
factors through $\sym_2W \otimes_K M_d$. Since $\wedge^2 W$
is the kernel of $W\otimes_K W\rTo \sym_2 W$, this is the same
as saying that the induced map $\wedge^2 W\otimes_K M_d\rTo M_{d+2}$
is 0, or again that the map 
$\phi^2:\ \Hom_K(E_v,M_d)\rTo \wedge^2 V\otimes_K \Hom_K(E_v,M_{d+2})$,
is zero. 
This last is equivalent to $\RR(M)$ being a complex.
As the whole construction is reversible, we are done.
\Box

As a first step in extending $\RR$ to all complexes, 
we consider the case of a module regarded as a complex with a single
term, but in arbitrary cohomological degree. Let $M$ be
an$S$-module, regarded as a complex concentrated in cohomological
degree 0. Then $M\langle a\rangle$ is a complex concentrated in 
cohomological degree $-a$, and we set 
$$
\RR(M\langle a\rangle)=\RR(M)\langle a\rangle.
$$
Now consider the general case of
a complex of graded $S$-modules
$$
\FF:\qquad \cdots\rTo F^i\rTo F^{i+1}\rTo\cdots .
$$
Applying $\RR$ to each $F^i$, regarded as a complex concentrated
in cohomological degree $i$, we get a double complex,
and we define $\RR(\FF)$ to be the total complex of this double complex.
Thus $\RR(\FF)$ is the total complex of 
$$
\diagram
&&\uTo && \uTo\\
\dots&\rTo & \Hom_K(E,(F^{i+1})_j)
                      &\rTo & \Hom_K(E,(F^{i+1})_{j+1})&\rTo&\dots\\
&&\uTo      &     &   \uTo&\\
\dots&\rTo&\Hom_K(E,(F^{i})_j)
                      &\rTo & \Hom_K(E,(F^{i})_{j+1})&\rTo&\dots\\
&&\uTo&&\uTo
\enddiagram,
$$
where the vertical maps are induced by the differential of $\FF$ and
the horizontal complexes are the complexes $\RR(F^i)$ defined above.
As $E$-modules we have 
$$
(\RR\FF)^k=\sum_{i-j=k} \Hom_K(E,(F^i)_j)
$$
where $(F^i)_j$ is regarded as a vector space concentrated
in degree $j$. Thus as a bigraded $E$-module, $\RR(\FF)=\Hom_K(E,\FF)$, and
the formula for the graded components is
$$
\RR(\FF)^i_j=\sum_m \Hom_K(E_{j-m} ,\FF^{i-m}_m).
$$

The functor $\RR$ has a left adjoint $\LL$ defined in an
analogous way by tensoring with $S$: on a graded $E$-module $P=\sum P_j$
the functor $\LL$ takes the value 
$$
\LL(P):\qquad \dots\rTo S\otimes_K P_j\rTo S \otimes_K P_{j-1}\rTo \dots ,
$$
where the map takes $s\otimes p$ to $\sum_ix_is\otimes e_ip$
and the term $S\otimes_K P_j$ has cohomological degree $-j$.
If $\GG$ is a complex of graded $E$-modules, then we can apply
$\LL$ to each term to get a double complex, and we define
$\LL(\GG)$ to be the total complex of this double complex,
so that 
$$
\LL(\GG)^k=\sum_{i-j=k}S\otimes_K(G^i)_j
\quad\hbox{and}\quad
\LL(\GG)^i_j=\sum_{m}S_{j-m}\otimes_K(G^{i+m})_m.
$$

To see that $\LL$ is the left adjoint of $\RR$ we proceed
as follows. First, if $M$ and $P$ are
left modules over $S$ and $E$ respectively, then
$$
\Hom_S(S\otimes_K P, M)=\Hom_K(P,M)=\Hom_E(P, \Hom_K(E,M)).
$$
If now $\MM$ and $\PP$ are complexes of graded modules over $S$ and
$E$, we must prove that $\Hom_S(\LL(\PP),\MM)\iso
\Hom_E(\PP,\RR(\MM))$, where on each side we take the maps of modules
that preserve the internal and cohomological degrees and commute with
the differentials.  As a bigraded module, $\LL(P)=S\otimes_KP$, and
similarly for $\RR$.  Direct computation shows that these maps of
complexes correspond to the maps of bigraded $K$-modules
$$
\phi=(\phi^i_j)\in \Hom_{\rm bigraded\ vector \ spaces}(\PP,\MM)
$$ 
such that $\phi^i_j:\ P^i_j\to M^{i-j}_j$ and
$$
\phi d-d\phi=(\sum_s x_s\otimes e_s)\phi,
$$
where $(\sum_s x_s\otimes e_s)\phi$ takes an element $p\in P^i_j$
to
$(-1)^i\sum_s x_s\phi(e_sp)$.
We have proved:

\theorem {BGG theorem}
{\bf (Bernstein-Gel'fand-Gel'fand [1978])} The functor $\LL$,
from the category of complexes of graded $E$-modules to the category
of complexes of graded $S$-modules, is a left adjoint to the functor
$\RR$.\Box

It is not hard to compute the homology of the complexes
produced by $\LL$ and $\RR$:

\proposition{koszul homology} If $M$ is a graded $S$-module and
$P$ is a graded $E$-module then
\item{$a)$} $\H^i(\RR(M))_j=\tor^S_{j-i}(K,M)_j$
\item{$b)$} $\H^i(\LL(P))_j=\ext_P^{j-i}(K,P)_j$

\proof  The $j-i^\th$ free module in the free resolution
of $K$ over $E$ is $(\Sym_{j-i}(W))^*\otimes_K E$, which is generated
by the vector space $(\Sym_{j-i}(W))^*$ of degree $i-j$. We can
use this to compute the right hand side of the equality in $b)$: the
$j^\th$ graded component of the module of homomorphisms of this into $P$
may be identified with $\Sym_{j-i}(W)\otimes_K P_{i}$. The differential
is the same as that of $\LL(P)$, and part $b)$ follows. Part $a)$ is
similar (and even more familiar, from Koszul cohomology.)\Box

It follows that the exactness of $\RR(M)$ or $\LL(P)$ are familiar
conditions. First the case a module over the symmetric algebra:

\corollary{exactness criterion} 
\item{$a)$} If $M$ is a finitely generated graded $S$-module,
then 
the truncated complex
$$
\RR(M)_{\geq d}: \Hom_K(E,M_d)\rTo \Hom_K(E,M_{d+1})\rTo\dots
$$
is acyclic (that is, has homology only at $\Hom_K(E,M_d)$) if
and only if $M$ is $d$-regular. 

\proof By \ref{koszul homology} applied to $M_{\geq d}$ the 
given sequence is acyclic if and only if $M_{\geq d}$ has
linear free resolution.\Box

Since any linear complex is of the form $\LL(P)$ for
a unique graded $E$-module $P$ it is 
perhaps most interesting to interpret part $b)$ of
\ref{koszul homology} as a statement about linear complexes
over $S$. The result below is implicitly used in Green's [1999] proof 
of the Linear Syzygy Conjecture. 

We call a right bounded linear complex 
$$
\GG:\quad \dots\rTo G^{-2}\rTo G^{-1}\rTo^{\phi} G^0
$$ 
{\it irredundant\/} if it is a subcomplex of the 
minimal free resolution of $\coker(\phi)$ (or equivalently of
any module whose presentation has linear part equal to $\phi$.)
(Eisenbud-Popescu [1999] called this property {\it linear exactness\/},
but to follow this usage would risk overusing the adjective
``linear''.) 

\corollary{linear exactness} Let $\GG$ be a minimal linear complex of
free $S$-modules ending on the right with $G_0$ as above,  and let
$P^*$ be the
$E$-module such that
$\LL(P^*)=\GG$. The complex $\GG$ is 
irredundant if and only the module $P$ is generated by $P_0$.
The complex $\GG$ is the linear part of a minimal free resolution
if and only if the module $P$ is linearly presented.

\proof 
Let $\phi: G^{-1}\rTo G^0$ be the differential of $\GG=\LL(P)^*$, let
$$
\FF:\quad \dots\rTo F^{-2}\rTo F^{-1}\rTo^{\phi} G^0
$$ 
be the minimal free resolution of $\coker(\phi)$,
and let $\kappa:\GG\rTo \FF$ be a comparison map lifting the 
identity on $G^0$. 
(This comparison map is unique
because $\FF$ is minimal and $\GG$ is linear.) By induction
one sees that the comparison map is an
injection if and only if $\H^{i}\GG_{-i}=0$ for all $i<0$, 
and it is an isomorphism onto the linear part of $\FF$ if and
only if in addition $\H^{i}\GG_{1-i}=0$ for all $i<0$. 
\ref{koszul homology} shows that the first condition is
satisfied if and only if $P^*$ injects into a direct sum
of copies of $E$, while both conditions are true if and only
if the minimal injective resolution begins with 
$$
0\rTo P^*\rTo E^a\rTo E(-1)^b
$$
for some numbers $a,b$. Dualizing, we get the desired conditions
on $P$.
\Box

We now return to the BGG-correspondence.
Both the functors $\LL$ and $\RR$ preserve
mapping cones and homotopies of maps of complexes. For mapping cone 
this is immediate. For the second note that two maps 
$f,g:\FF\to \GG$
of complexes are homotopic if and only if the induced map
from $\FF$ to the mapping cone of
$f-g$ is split. This condition is preserved by any additive functor
that preserves mapping cones. 

Recall that a free resolution of a right bounded complex 
$$
\MM:\qquad \dots\rTo M^{i-1}\rTo M^{i}\rTo M^{i+1}\rTo \dots
$$
of graded $S$-modules is a graded free complex $\FF$ with a
morphism $\FF\rTo \MM$, homogeneous of degree 0, which 
induces an isomorphism on homology.
We say that $\FF$
is {\it minimal\/} if $K\otimes_S\FF$
has trivial differential. Every right bounded complex
$\MM$ 
of finitely generated modules has a minimal free resolution,
unique up to isomorphism. It is the minimal part of any free resolution.

The functors $\LL$ and $\RR$ give a general construction of resolutions. 

\theorem{ecoh-thm1} 
For any complex of graded $S$-modules $\MM$, the complex
$
\LL\RR(\MM)
$
is a free resolution of $\MM$ which surjects onto $\MM$; and
for any complex of graded $E$-modules $\GG$,
the complex
$
\RR\LL(\GG) 
$
is an injective resolution of $\GG$ into which $\GG$ injects.

An immediate consequence is:

\corollary{BGG-result} The functors $\RR$ and $\LL$ 
define an equivalence $D^b(S\hbox{-mod})\iso D^b(E\hbox{-mod})$.

\proof The derived category $D^b(S\hbox{-mod})$ of bounded complexes
of finitely generated $S$-modules is equivalent to the derived
category of complexes of finitely generated $S$-modules with
bounded cohomology (that is, having just finitely many 
cohomology modules), see for example Hartshorne [1977], III Lemma 12.3,
and similarly for $E$. By
\ref{ecoh-thm1}, the functors $\LL$ and $\RR$ carry
bounded complexes into complexes with bounded cohomology,
and $\LL\RR, \RR\LL$ are both equivalent to the identity.
\Box

\noindent {\sl Proof of \ref{ecoh-thm1}.\/} 
The proofs of the two statements are similar,
so we treat only the first. (A slight simplification
is possible in the second case since
finitely generated modules over $E$ have
finite composition series.)

Because $\LL$ is the left adjoint functor of $\RR$ there is a 
natural map $\LL\RR(\MM)\rTo \MM$ adjoint to the identity map
$\RR(\MM)\rTo \RR(\MM)$. We claim that this is a surjective
quasi-isomorphism. 

To see that it is a surjection, 
consider a map $\phi:\; \MM\rTo \MM'$ such that the composite
$\LL\RR(\MM)\rTo \MM\rTo \MM'$ is zero. It follows that
the adjoint composition $\RR(\MM)\rTo \RR(\MM)\rTo \RR(\MM')$
is also zero, and since the first map is the identity,
we get $\RR(\phi)=0$. Since $\RR$ is a faithful functor, 
$\phi=0$, proving surjectivity.

The functor $\LL$ preserves direct limits because it is a left 
adjoint, while the functor $\RR$ preserves direct limits because
$E$ is a finite dimensional vector space. Thus it suffices
to prove our claim in the case where $\MM$ is a bounded complex
of finitely generated $S$-modules. 

If $\MM$ has the form 
$$
\MM:\qquad \dots\rTo M^d\rTo 0\rTo \dots
$$
then $\MM$ admits $M^d\langle -d\rangle$ 
(that is, the module $M^d$ considered as
a complex concentrated in cohomological degree $d$) as a subcomplex, and
the quotient is a complex of smaller length. Using the ``five lemma''
the naturality of the map $\LL\RR(\MM) \rTo \MM$, and the exactness
of the functor $\LL\RR$, the claims
will follow, by induction on the length of the complex, 
from the case
where $\MM$ has the form $M\langle -d\rangle$ for some finitely generated
graded $S$-module $M$ and integer $d$. This reduces immediately
to the case $d=0$. 

It thus suffices to to see that 
$\LL\RR(M)\rTo M$ is a quasi-isomorphism when $M$ is a finitely
generated graded $S$-module. Now $\RR(M)$ is the linear complex
$\Hom_K(E,M_0)\to\Hom_K(E,M_1)\to\cdots$, so
$\LL\RR(M)$ is the total complex
of the following double complex:
$$\diagram
     &       & \uTo                   &    \\
\cdots& \rTo  &S\otimes_K\Hom_K(K,M_1) &\rTo& 0                      \\
     &       & \uTo                   &    & \uTo \\
\cdots& \rTo  &S\otimes_K\Hom_K(V,M_0) &\rTo& S\otimes_K\Hom_K(K,M_0)&\rTo& 0.
\enddiagram
$$
In this picture the terms below what is shown are all zero.
The terms of cohomological degree 0 in the total complex
are those along the diagonal going northwest from 
$S\otimes_K\Hom_K(K,M_0)$. The generators of 
$S\otimes_K\Hom_K(K,M_0)$ have internal degree 0, while those
of $S\otimes_K\Hom_K(K,M_1)$ have internal degree 1, etc.

The $d^\th$ row of this double complex is $S\otimes_K\Hom_K(E,M_d)$,
which is equal to the complex obtained by tensoring the 
Koszul complex 
$$
\dots \to S\otimes_K\wedge^2 W \to S\otimes_K W \to S \to 0
$$
with $M_d$. It is thus acyclic, its one cohomology module being
$M_d$, in cohomological degree $0$. The spectral sequence 
starting with the horizontal cohomology of the double complex 
thus degenerates, and we see that the cohomology of the total
complex $\LL\RR(M)$ is a graded module with
component of internal degree equal to $M_d$, 
concentrated in cohomological degree 0.
Thus
$\LL\RR(M)$ is acyclic and the Hilbert
function of $\H^0(\LL\RR(M))$ is the same as that of $M$.
As $\LL\RR(M)$ has no terms in
positive cohomological degree, and $M$ is in cohomological degree 0,
the surjection $\LL\RR(M)\rTo M$ induces a surjection
$\H^0(\LL\RR(M))\rTo M$, and we are done. (One can show
that $\LL\RR(M)$ is the tensor product, over $K$,
of the Koszul complex and $M$, the action of $S$ being
the diagonal action, but the isomorphism is complicated
to write down.) \Box

Though the statement of \ref{ecoh-thm1} has an attractive simplicity,
it is not very useful in this form because the resolutions that are
produced are highly nonminimal (for example the free resolutions
produced over $S$ are nearly always infinite). 
\ref{reciprocity} shows that a modification of this construction gives at
least an important part of the minimal free resolution.

\section{linear part section} The Linear Part of a Complex

If $A$ is a matrix over $E$ then we define the {\it linear part\/},
written $\lin(A)$, to be the matrix obtained by erasing all the terms
of entries of $A$ that are of degree $>1$.  For example, if
$a,b,c,d$ are linear forms of $E$, then 
the linear part of
$$\pmatrix{
a&0\cr
bc&d}
\qquad\hbox{is}\qquad
\pmatrix{
a& 0\cr
0&d}.
$$
Taking the linear part is a functorial operation
on maps (see \ref{linear part and tor} below), but
taking the linear part of a matrix does not
always commute change of basis. 
For example, if $a,b,c$ are linear
forms,
$$ d=\pmatrix{ 
a & 0 \cr
0 & b
},
\qquad \hbox{and}\qquad
 e=\pmatrix{1 & c \cr
               0 & 1}, $$
then
$\lin(de) \neq \lin(d)e. $

Suppose that $e:G\rTo H$ is a second map of free modules and that the
composition $ed=0$. It {\it need not\/} be the case that that
$\lin(e)\lin(d)=0$; but if we assume in addition that $d(F)$ is in the
maximal ideal times $G$ and $e(G)$ is in the maximal ideal times $H$,
then $\lin(e)\lin(d)=0$ does follow. Thus if $\FF$ is a
minimal free complex over $E$ we can define a new complex $\lin(\FF)$
by replacing each differential $d$ of $\FF$ by its linear part,
$\lin(d)$. Note that $\lin(\FF)$ is the direct sum of complexes
$\FF^{(d)}$ whose $i^\th$ term $\FF^{(d)}_i$ is a direct sum of copies
of $E(d+i)$ and whose maps are of degree 1. In general, we define
the linear part of any free complex $\FF$ to be the linear part of
a minimal complex homotopic to $\FF$. 

\theorem{linear dominance} Let $\FF$  be a free or injective resolution 
of a finitely generated module over the exterior algebra $E$. The
linear part of $\FF$ is eventually exact.

\proof We treat only the case where $\FF$ is an injective resolution;
by duality, the statement for a free resolution is equivalent.
By \ref{linear part and tor} the linear part of $\FF$ is
the value of $\RR$ on the $S$-module $\Ext_E^\bullet(K,M)$. 
Since any finitely generated $S$ module has finite regularity
(see Eisenbud-Goto [1984]), it suffices by \ref{exactness criterion} 
to show that 
$\Ext_E^\bullet(K,M)$ is a finitely generated
$S$-module. This was done by Aramova, Avramov,
and Herzog [2000]. For the reader's convenience we repeat the
argument: we prove that $\Ext_E^\bullet(K,M)$ is a finitely generated
$S$-module by induction on the length
of $M$. If $M=K$,
then $\Ext_E^\bullet(K,M)$ is free of rank 1 over $S$. If $M'$ is a
proper submodule of $M$ then from the exact sequence
$$
0\rTo M' \rTo M \rTo M/M' \rTo 0
$$
we get an exact triangle of $S$-modules
$$
\diagram
\Ext_E^\bullet(K,M/M') & \rTo & \Ext_E^\bullet(K,M')\\
\hskip 6em\luTo &\hskip 6em \ldTo\\
&\Ext_E^\bullet(K,M) &\\
\enddiagram.
$$
The two $S$-modules in the top row are finitely generated by
induction, and thus $\Ext_E^\bullet(K,M)$ is finitely generated too.
\Box

If $M$ is an $E$-module, then we write $\lin(M)$ for the cokernel of
$\lin(d)$, where $d$ is the map in a minimal free presentation of
$M$. We can further define a family of modules connecting $M$ and
$\lin(M)$ as follows: Let $d$ be a minimal free presentation of $M$,
choose a representation of $d$ as a matrix, and let $e_1,\dots,e_v$ be
a basis of $V$.  Let $d'$ be the result of substituting $te_i$ for
$e_i$ in the entries of $d$, and then dividing each entry by $t$. The
entries of $d$ have no constant terms because $d$ is minimal, and it
follows that $d'$ is a matrix over $K[t]\otimes_K E$. Let $M'$ be
the cokernel of $d'$. It has
fibers $M$ at $t\neq 0$ and $\lin(M)$ at 0.
The module $M'$ may not be flat over $K[t]$,
but the module $M'[t^{-1}]$ is flat over $K[t,t^{-1}]$:
in fact, it is 
isomorphic to the module obtained
from the trivial family $K[t,t^{-1}]\otimes_K M$ by pulling back 
along the automorphism $e_i\mapsto te_ic$ of $E$.

\corollary{deformations} If $M$ is a finitely generated $E$ module,
then any sufficiently high syzygy $N$ of $M$ is a flat deformation
of its linear part $\lin(N)$.

\proof If $N$ is a sufficiently high syzygy, then by
\ref{linear dominance} the linear part of the minimal resolution of $N$ 
is the resolution of $\lin(N)$, so that
(with the notation of the preceding paragraph) this free resolution of
$N$ lifts to a free resolution of $\lin(N')$ over $K[t]\otimes_K E$.
Thus $N'$ is flat, and the result follows.
\Box

\example{eg1} 
It is sometimes not so obvious what the linear part of the minimal
version of a complex will be, and in particular it may be hard
to read from the linear terms in a nonminimal version. For example,
suppose that $W$  has dimension 2 and that $x,y\in W$ is
a dual basis to $a,b\in V$. Consider the complex
$$
\GG:\quad 0\rTo S/(x,y^2)\rTo^x S/(x^2,y)(1)\rTo 0
$$
where the notation means that the class of 1 goes to the class of $x$.

Applying $\RR$ to $\GG$, we get the double complex
$$
\diagram
0& \rTo  &   E(1)  & \rTo^a & E       &\rTo          &0    \cr   
 &       &   \uTo   &        &\uTo^1   &            &\uTo              \cr
 &       &     0    &  \rTo  &  E      &\rTo_b      &E(-1)& \rTo  &   0\cr
\enddiagram
$$
whose total complex is 
$$
\RR(\GG)=\FF:\quad 0\rTo E(1)\oplus E
\rTo^{\pmatrix{
a&1\cr
0&b
}}
E\oplus E(-1)\rTo 0.
$$
Despite the presence of the linear terms
in the differential of $\FF$, the minimal complex 
$\FF'$ homotopic
to $\FF$ is 
$$
\FF':\quad 0\rTo E(1)\rTo^{ab} E(-1)\rTo 0
$$
so the differential of $\lin(\FF)$ is 0.

Fortunately, we can construct the linear part of a complex directly and
conceptually, without passing to a minimal complex or to matrices.
First note that if $\GG$ is a minimal free complex, 
then giving its linear part is equivalent, by \ref{basic correspondence},
to giving the maps
$\phi_i:K\otimes_E G^i\rTo V\otimes_K K\otimes_E G^{i+1}$ 
corresponding to the
linear terms in the differential of $\GG$.
If $\FF$ is any free complex homotopic to $\GG$,
then $K\otimes_E G^i=\H^i(K\otimes\FF)$. We will construct
natural maps
$\psi_i:\H^i(K\otimes\FF)\rTo V\otimes\H^{i+1}(K\otimes_E\FF)$, 
and prove that $\psi_i=\phi_i$.

We identify $S$ with $\Ext_E(K,K)$ and use the well-known 
$\Ext_E(K,K)$-module structure on $\H(K\otimes_E \FF)$. 
To formulate this explicitly,
we make use of the exact sequence 
$$
\eta:\qquad 0\rTo V\rTo E/(V)^2\rTo K\rTo 0.
$$
The extension class 
$$
\eta\in\Ext^1_E(K,V)=\Ext^1_E(K,K)\otimes_K V=Hom_K(W,\Ext^1_E(K,K))
$$
corresponds to the inclusion $W=\Sym_1W\subset \Sym W$.
Since $\FF$ is a free complex, the sequence
$\eta\otimes_E\FF$ is an exact sequence of complexes, and we obtain
the homomorphism $\psi_i$ as the connecting homomorphism
$$\H_i(K\otimes_E \FF)\rTo H_{i+1}(V\otimes_E \FF)=
V\otimes_K\H_{i+1}(K\otimes_E\FF).
$$

\theorem{linear part and tor} If $\FF$ is a complex of free
modules over $E$,
then 
$$
\lin(\FF)=\RR(\H^\bullet(K\otimes_E \FF)),
$$
where the $S$-module structure
on
$\H^\bullet(K\otimes_E \FF)$ is given by the action of $\Ext_E(K,K)$.

\proof From the definition of $\psi_i$ we see that it is
preserved by homotopy equivalences of free complexes, so we may
assume that $\FF$ is minimal. 
To simplify the notation, we set
$\bar F_i=K\otimes_E F_i$, so that $F_i= E\otimes_K\bar F_i$.
Suppose $a\in \bar F_i$ goes to $b\in V\otimes_K \bar F_{i+1}$ under $\phi_i$,
the linear part of the differential $d$ of $\FF$. To compute $\psi_i(a)$ we
first lift $a$ to the element 
$$
1\otimes a\in E/(V)^2\otimes_E F_i
$$
and then apply $E/(V)^2\otimes_E d$. The result is $b$, regarded
as an element of $E/(V)^2\otimes_E F_i$ via the inclusion
$V\subset E/(V)^2$. Since this is the inclusion that defines
$\eta$, we are done.\Box

To understand the linear parts of complexes obtained from
the functor $\RR$, we will
employ a general result:
if the vertical differential of suitable double complex splits,
then the associated total complex is homotopic to one built from
the homology of the vertical differential in a simple way.

\lemma{degenerate double complex} Let $\FF$ be a double complex 
$$
\diagram
&&\uTo && \uTo\\
\dots&\rTo & F^{i+1}_j
                      &\rTo^{d_{hor}} & F^{i+1}_{j+1}&\rTo&\dots\\
&&\uTo^{d_{vert}}      &     &   \uTo_{d_{vert}}&\\
\dots&\rTo&F^{i}_j
                      &\rTo_{d_{hor}} & F^{i}_{j+1}&\rTo&\dots\\
&&\uTo&&\uTo
\enddiagram,
$$
in some abelian category such that $F^i_j=0$ for $i\ll 0$.
Suppose that the vertical differential of $\FF$ splits, so that for 
each $i,j$
there is a decomposition
$F^i_j=G^i_j\oplus d_{vert}G^{i-1}_j\oplus H^i_j$
such that the kernel of $d_{vert}$ in $F^i_j$ is
$H^i_j\oplus d_{vert}G^{i-1}_j$,
and such that $d_{vert}$ maps $G^{i-1}_j$ isomorphically to 
$d_{vert}(G^i_j)$. 
If we write $\sigma:\; F^i_j\to H^i_j$ for the projection corresponding
to this decomposition, and 
$\pi:\; F^i_j\to d_{vert}G^{i-1}_j\to G^{i-1}_j$ 
for the composition of the projection with the inverse of
$d_{vert}$ restricted to $G^{i-1}_j$, then the total complex of $\FF$ is
homotopic to the complex
$$
\dots\rTo \oplus_{i+j=k}H^i_j\rTo^d \oplus_{i+j=k+1}H^i_j\rTo\dots
$$
with differential 
$$
d=\sum_{\ell\geq 0}\sigma (d_{hor} \pi)^\ell d_{hor} 
$$

\proof We write $d_{tot}=d_{vert}\pm d_{hor}$ for the differential
of the total complex. Note first that $\sigma (d_{hor} \pi)^j d_{hor}$ takes
$H^i_j$ to $H^{i-\ell}_{j+1+\ell}$. Since $F^{i-\ell}_{j+1+\ell}=0$ 
for $\ell>>0$,
the sum defining $d$ is finite.

Let $F$ denote $\FF$ without the differential, that is, as a bigraded module. 
We will first show that 
$F$ is the direct sum of the three components 
$$
G=\oplus_{i,j} G^i_j,
\quad 
d_{tot}G, 
\quad\hbox{and   } 
H=\oplus_{i,j} H^i_j
$$
and that $d_{tot}$ is a monomorphism on $G$.

The same statements, with $d_{tot}$ replaced by
$d_{vert}$, are true by hypothesis. 
In particular, any element of $F$ is a sum of elements of the form
$g'+d_{vert}g+h$ with $g'\in G^i_j,\ g\in G^{i-1}_j$ and $h\in H^i_j$
for some $i,j$. Modulo $G+d_{tot}G+H$ this element can be written as
$d_{hor}(g)\in F^{i-1}_{j+1}$. As $F^s_t=0$ for $s<<0$, we may
do induction on $i$, and
assume that $d_{hor}g \in G+d_{tot}G+H$, so we see
that $F=G+d_{tot}G+H$.

Suppose 
$$
g'\in G=\oplus_{i+j=\ell} G^i_j,
\quad 
g\in G=\oplus_{i+j=\ell-1} G^i_j,
\quad\hbox{and   } 
h\in H=\oplus_{i+j=\ell} H^i_j
$$
and $g'+d_{tot}g+h=0$; we must show that $g=g'=h=0$.
Write $g=\sum_{k=a}^b  g^{k-1}_{\ell-k}$ with
$g^s_t\in G^s_t$. 
If $b-a=-1$ then
$d_{tot}=0$ and the desired result is
a special case of the hypothesis. In any case, there is no
component of $g$ in $G^{b}_{\ell-b-1}$ so
the
component of $d_{tot}g$ in
$G^{b}_{\ell-b}$ is equal to $d_{vert}g^{b-1}_{\ell-b}$. From the
hypothesis we see that $d_{vert}g^{b-1}_{\ell-b}=0$, 
so $g^{b-1}_{\ell-b}=0$, and we are done by induction
on $b-a$. This shows that $F=G\oplus d_{tot}G\oplus H$ and that
$d_{tot}$ is an isomorphism from $G$ to $d_{tot}G$.

The modules $G\oplus d_{tot}G$ form a 
double complex contained in $\FF$ that we will call $\GG$.
Since $d_{tot}:G\rTo d_{tot}G$ is
an isomorphism, the total complex $tot(\GG)$ is split exact.
It follows that the total complex $tot(\FF)$ is
homotopic to $tot(\FF)/tot(\GG)$, and the modules of
this last complex are isomorphic to $\oplus_{i+j=k}H^i_j$. 
We will complete
the proof by showing that the induced differential on
$tot(\FF)/tot(\GG)$ is the differential $d$ defined above.

Choose $h\in H^i_j$. The image of $h$ under the induced 
differential is the unique element $h'\in H$ such that
$d_{tot}h\equiv h'\ (\hbox{mod }G+dG)$.
Now
$$
d_{tot}h
=d_{hor}h
\equiv \sigma d_{hor} h + (d_{vert}\pi) d_{hor} h\ (\hbox{mod }G).
$$
However, 
$$
d_{vert}\pi
\equiv  d_{hor}\pi 
\equiv \sigma (d_{hor}\pi) + d_{vert}\pi (d_{hor}\pi)
\ (\hbox{mod }G+d_{tot}G)
$$
Continuing this way, and using again the fact that $F^i_j=0$ for
$i<<0$ we obtain 
$$
d_{tot}h
\equiv \sum_\ell\sigma (d_{hor}\pi)^\ell d_{hor}h 
\ (\hbox{mod }G+d_{tot}G)
$$
as required.
\Box

We apply \ref{linear part and tor} to complexes of the form
$\RR(\FF)$:

\corollary{linear part 1} If $\GG$ is a left-bounded 
complex of graded $S$-modules,
then 
$$
\lin(\RR(\GG))=\oplus_i\RR(\H^i(\GG)),
$$
where $\H^i(\GG)$ is regarded as a complex of one term, concentrated
in cohomological degree $i$. A similar statement holds for
the linear part of $\LL(\GG)$ when $\GG$ is a left-bounded 
complex of graded $E$-modules.

\proof As $\GG$ is a left-bounded complex of finitely generated
modules, the double complex $(1)$ whose total complex is $\RR(\GG)$
satisfies the conditions of \ref{degenerate double complex}. The 
bigraded module underlying
$\RR(\H^i(\GG))$ is precisely the module $H$ of 
\ref{degenerate double complex},
and the differential is the map $\sigma d_{hor}$ restricted to $H$.
This is a linear map. But the other terms in the sum
$d=\sum_\ell \sigma(\pi d_{hor})^\ell d_{hor}$ all involve two or
more iterations of $d_{hor}$, and are thus represented by matrices 
whose entries have degree at least 2.
\Box

\medskip
\noindent{\bf\ref{eg1} Continued:} 
Note that the homology of $\GG$ is 
$\H^\bullet(\GG)=K(-1)\oplus K(1)\langle-1\rangle$.
Since the two copies of $K$ have degree differing by more
than 1, the $\Ext_E^\bullet(K,K)$-module structure is trivial.
We may write
$\lin(\FF')=\RR(K(-1))\oplus\RR(K(1))\langle -1\rangle$
as required by \ref{linear part 1}.

Here is the promised information about the minimal resolution of a module:

\theorem{reciprocity} $a)$ {\bf Reciprocity}: If $M$ is a finitely
generated graded $S$-module and $P$ is a finitely generated graded
$E$-module, then $\LL(P)$ is a free resolution of $M$ 
if and only if $\RR(M)$ is an injective resolution of $P$. \hfill \break
\indent $ b)$ More generally, for any minimal
bounded complex of finitely generated 
graded $S$-modules $\MM$, the linear part of the minimal
free resolution of $\MM$ is
$
\LL(\H^\bullet(\RR(\MM));
$
and
for any minimal bounded complex of finitely generated 
graded $E$-modules $\GG$,
the linear part of the minimal injective resolution of $\GG$
is
$
\RR(\H^\bullet\LL(\GG)).
$

\proof 
The two parts of $b)$ being
similar, we prove only the first statement. By \ref{ecoh-thm1}
the complex $\LL\RR(\FF)$ is a free resolution. The complex
$\RR(\FF)$ is left-bounded because $\FF$ is bounded and contains
only finitely generated modules. Thus we may apply 
\ref{linear part 1}, proving the first statement.

For the reciprocity statement $a)$, suppose that $\LL(P)$ is a minimal free
resolution of
$M$. By part b) the linear part of the minimal 
injective resolution of $P$ is $\RR(\H^\bullet(\LL(P)))$. Since
$\LL(P)$ is a resolution of $M$, this is $\RR(M)$. All the 
terms of cohomological degree $d$ of this complex have degree $-d$,
so there is no room for nonlinear differentials, and the linear part
of the resolution is the resolution. \Box

\section{tate} Sheaf cohomology and exterior syzygies

In this section we establish a formula for the free
modules that appear in resolutions over $E$. 
Because $E$ is Gorenstein, it is 
natural to work with doubly infinite resolutions:

A {\it Tate resolution\/} over $E$ is a doubly infinite free
complex 
$$
T:\qquad \dots \rTo T^d\rTo T^{d+1}\rTo \dots
$$
that is everywhere exact. 

There is a Tate resolution naturally associated to a
coherent sheaf $\F$
on $\P(W)$, defined as follows. Let
$M$ be a finitely generated graded
$S$-module representing $\F$, 
for example $M=\oplus_{\nu \geq 0}\H^0(\F(\nu))$.
If $d\geq\reg(M)$, then by \ref{exactness criterion}
the complex $\RR (M_{\geq d})$ is acyclic. Thus
if $d>\reg(M)$ then, since $\RR(M_{\geq d})$ is minimal,
$\Hom_K(E,M_d)$ minimally covers
the kernel of the map $\Hom_K(E,M_{d+1})\rTo \Hom_K(E,M_{d+2})$

Fixing $d>\reg(M)$, we may complete $\RR(M_{\geq d})$ to
a minimal
Tate resolution $\TT(\F)$ by adjoining a free resolution of 
$$
\ker\bigl[\Hom_K(E,M_d)\rTo \Hom_K(E,M_{d+1})\bigr].
$$
Since any two modules representing $\F$ are equal
in large degree, the Tate resolution
is independent of which $M$ and which large $d$ is chosen,
and depends only on 
the coherent sheaf $\F$. It has the form
$$\eqalign{
\TT(\F)&: \qquad\cdots \to \cr
&T^{d-2}\to T^{d-1}\to \Hom_K(E, \H^0(\F(d)))\to
\Hom_K(E, \H^0(\F(d+1)))\cr
&\qquad\qquad\qquad\qquad\qquad\qquad\qquad\qquad\qquad\qquad\qquad\qquad
\to\cdots}
$$
where the $T_i$ are graded free $E$-modules.

The main theorem of this section expresses 
the linear part of this Tate resolution in terms of the 
$S$-modules $\oplus_e \H^j(\F(e))$ given by 
the (Zariski) cohomology of $\F$. We regard
$\oplus_e\H^j(\F(e))$ as an $S$-module concentrated in
degree $j$.

\theorem{sheaf cohomology} If $\F$ is a coherent sheaf
on $\P(W)$, then the
linear part of the Tate resolution $\TT(\F)$ is 
$\oplus_j \RR(\oplus_e \H^j(\F(e)))$.
In particular,
$$
T^e=\oplus_j \Hom_K(E,\H^j(\F(e-j))),
$$
where $\H^j(\F(e-j))$ is regarded as a vector space of 
internal degree $e-j$.

A special case of the theorem appears without proof as Remark 3 after
Theorem 2 in Bernstein-Gel'fand-Gel'fand [1978].  
The proof below could be extended to
cover the case of a bounded complex of coherent sheaves, replacing the
cohomology in the formula with hypercohomology. 
We will postpone the proof of \ref{sheaf cohomology}
until the end of this section. 

Rewriting the indices in \ref{sheaf cohomology}, we emphasize
the fact that we can compute any part of the cohomology of $\F$
from the Tate resolution. 

\corollary{comp of shf coho} For
all $j,\ell\in \Z$,
$$
\H^j(\F(\ell)) = \Hom_E(K,T^{j+\ell})_{-\ell}.
$$
\Box

\ref{comp of shf coho}  provides the basis for 
an algorithm
computing the cohomology of $\F$ with
any computer program that can provide free resolutions
of modules over the symmetric and exterior algebras, 
such as the program Macaulay2 of Grayson and Stillman
[{\tt http://www.math.uiuc.edu/Macaulay2/}].
For an explanation of the algorithm in practical terms, see
((**** Reference to the Macaulay2 book to come [****]))

To prove \ref{sheaf cohomology} we will use 
the reciprocity result \ref{reciprocity}. We actually prove 
a slightly more general version,
involving local cohomology. We write $\gm$ for the homogeneous maximal
ideal $SW$ of $S$, and for any graded $S$-module $M$ we write
$\H^j_\gm(M)$ for the $j^\th$ local cohomology module of $M$, regarded
as a graded $S$-module.

\theorem{local cohomology} Let $M$ be a graded $S$-module
generated in degree $d$, and having
linear free resolution $\LL(P)$. Let
$\FF:\ \cdots \to F^{-1}\to F^0$ be the minimal free resolution of $P$. 
The linear part
of $\FF$ is 
$$
\lin(\FF)=\oplus_j\RR(\H^j_\gm(M)),
$$ 
where $\H^j_\gm(M)$ is 
regarded as a complex with one term, concentrated in
cohomological degree $j$. In particular, we have 
$$
F^{-i}=\oplus_j\Hom_K(E, \H^j_\gm(M)_{-j-i}).
$$

\noindent{\sl Proof of \ref{local cohomology}.\/}
We compute the linear part of the free resolution of $P$ by 
taking the dual (into $K$) of the linear part of the injective 
resolution of $P^*$.
By \ref{reciprocity} the linear part of the injective resolution
of $P^*$ is $\RR(H^\bullet(\LL(P^*)))$. 
It follows at once from
the definitions that $\LL(P^*)=\Hom_S(\LL(P), S)$.
By once more \ref{reciprocity} $\LL(P)$ is the minimal free resolution
of $M$, so
$\H^\bullet(\LL(P^*))=\Ext_S^\bullet(M,S)$.
 Thus
the linear part of the free resolution of $P$ is
$[\RR\Ext_S^\bullet(M,S)]^*$, where $\Ext^j(M,S)$ is 
thought of as a module concentrated in cohomological degree
$j$.

Because $E^*=\omega_E=E\otimes \wedge^vW$ we have, for 
any graded vector space $D$, natural
identifications
$$\eqalign{
(\Hom_K(E,D))^*
&=(E^*\otimes_KD)^*\cr
&=E^*\otimes_K \wedge^vW^* \otimes_KD^*\cr
&=\Hom_K(E,D^*))\otimes_K\wedge^vV}
$$
(Here all the duals of $E$-modules are Hom into $K$.) If $D$
has the structure of a graded $S$-module then 
$D^*$ is again a graded $S$-module, and this becomes this is an 
isomorphism of graded $S$-modules. If we think of 
$D$ as a complex with just one term, in cohomological degree
$d$, then 
$\RR(D)^*=\RR(D^*\otimes_K\wedge^vV)$ where, to make
all the indices come out right, we must think of 
$D^*\otimes_K \wedge^vV=(D\otimes_K\wedge^vW)^*$ as a complex of one
term concentrated in cohomological degree $v-d$.

If we take $D=\Ext_K^\ell(M,S)$ then by local duality
$$\eqalign{
D^*&=(\Ext_K^\ell(M,S)\otimes_K\wedge^vW\otimes_K \wedge^vV)^*\cr
&=(\Ext_K^\ell(M,\omega_S)\otimes_K \wedge^vV)^*\cr
&=\H^{v-\ell}_\gm(M)\otimes \wedge^vW.}
$$
Thus
$$\eqalign{
\RR(\Ext_K^\ell(M,S)^*)&=\RR(\H^{v-\ell}_\gm(M)
                         \otimes_K\wedge^vW)\otimes\wedge^vV\cr
&=\RR(\H^{v-\ell}_\gm(M)).
}
$$
where $\H^j_\gm(M)$
is regarded as a complex with just one term, of
cohomological degree $-j$, as required.\Box

\noindent{\sl Proof of \ref{sheaf cohomology}.\/}
For each $i=0, \dots, v-1$ we write $\H^i$ for the
cohomology module $\oplus_{d=-\infty}^\infty\H^i(\F(d))$.
If we choose $d\geq\reg(\H^0_{\geq 0})$  as in the
construction of $\TT(\F)$, the module $M:=\H^0_{\geq d}$
has a linear free resolution, so we may 
apply \ref{local cohomology}. We deduce that
the linear part of the  free resolution of 
$P:=\ker [\Hom_K(E,\H^0(\F(d)))\rTo \Hom_K(\H^0(\F(d+1)))]$
is 
$
\lin(\FF)=\oplus_j\RR(\H^j_\gm(M)).
$
If we insist that $d>\reg(\H^0_{\geq 0})$
then $\H^0_\gm(M)=0$. From the exactness of the
sequence
$$
0\rTo \H^0_\gm(M)\rTo 
M\rTo
\oplus_{d=-\infty}^\infty \H^0(\F(d))\rTo
\H^1_\gm(M)\rTo 0
$$
it follows that
the local cohomology module $\H^{1}_\gm(M))$
agrees with the global cohomology module
$\H^0$ in all degrees strictly less than d,
and of course we have
$\H^i=\H^{i+1}_\gm(M)$. This concludes the proof.
\Box

As a further application of \ref{reciprocity} we have:

\corollary{lin cplx} If $\F$ is a coherent sheaf
on $\P(W)$ with
Tate resolution 
$$
\TT(\F):\qquad \cdots\rTo T^{-1}\rTo T^0\rTo T^1\rTo \cdots
$$
then the linear free complex of $S$-modules
$\GG=\LL(\ker(T^0\rTo T^1))$
satisfies 
$$
\H^i(\GG)=\oplus_{j\geq -i}\H^i(\F(j));\quad 
\H^{n-i}\Hom_S(\GG,S)=\oplus_{j< -i}\H^i(\F(j))^*\otimes \wedge^vV
$$
as $S$-modules.
In particular
the sheafified complex $\tilde \GG$ is a linear monad for $\F$ with at most
$2n+1$ terms.

\proof Set $P:= \ker(T^0\rTo T^1)$. By part $b)$ of \ref{reciprocity} the linear part
of the injective
resolution $T^0 \rTo T^1 \rTo \ldots$ of 
$P$ is the sum of the
linear complexes
$\RR(\H^i(\GG))$. But in \ref{sheaf cohomology}
we identified the
linear part of this resolution as 
$\RR(\oplus_{j\geq -i}\H^i(\F(j))^*).$ 
\ref{basic correspondence} and a comparison
of indices complete the
proof of the first formula.

The proof of the second formula follows similarly once one
observes that $\Hom(\GG, S)=\LL(\Hom_K(P,K))$; that 
the injective resolution of $\Hom_K(P,K)$ is 
$\Hom(\TT(\F),K)\langle 1\rangle_{\geq 0}$; and that 
the terms with $\H^i$ on the right hand side of the desired equality
correspond to $(v-i)^\th$ linear strand of $\Hom(\TT(\F),K)$.\Box

\remark{} \ref{reciprocity} also implies a closely related 
interpretation of the dual complex of $\TT(\F)$: it is the
Tate resolution of the dual of $\F$ in the derived 
category:
$$\Hom_K(\TT(\F),\wedge^v W) \cong 
\TT({\rm R}{\cal H}om(\F,\omega_\P\langle -n\rangle)).$$
This with the hypercohomology version of \ref{sheaf cohomology}
implies the usual statement of Grothendieck duality for
sheaves on $\P(W)$.

\example{elliptic quartic} Let $C$ be an elliptic quartic
curve in $\P^3$, and consider $\O_C$ as
a sheaf on $\P^3$. Write $\omega=\omega_E$ for
$\wedge^vW\otimes E(-v)$ as usual.
Computing cohomology 
one sees that $\TT(\O_C)$ has the form
$$
\cdots\rTo
\omega^8(2) 
\rTo 
\omega \oplus \omega^4(1)
\rTo^d
\omega^4(-1)\oplus\omega 
\rTo
\omega^8(-2)
\rTo \cdots
$$
The complex $\GG$ in the theorem above is obtained
by applying $\LL$ to the
kernel of the map marked $d$.

>From the given form
of the resolution, one easily computes the Hilbert function
of this kernel, and it follows that $\GG$ has the form
$$
0\rTo S^8(-2)\rTo S^{20}(-1) \rTo S^{16} \rTo S^4(1)\rTo 0,
$$
with the term $S^{16}$ in cohomological degree 0.

If $C \subset \P^3$ is taken to be Heisenberg invariant, say
$C = \{ x_0^2+x_2^2+\lambda x_1x_3=x_1^2+x_3^2+\lambda x_0x_2 = 0 \}$
for some $\lambda \in \AA^1_k$, then $d$ can be represented
by the matrix
$$
\pmatrix{ 0 &e_0 &         e_1   &     e_2   &    e_3         \cr
 e_0 &-\lambda e_1e_3  & e_2e_3  & 0 & e_1e_2+{\lambda^2 \over 2}e_0e_3 \cr
 e_1  &e_2e_3 & \lambda e_0e_2    &  -e_0e_3-{\lambda^2 \over 2} e_1e_2 & 0 \cr
 e_2 &0       & -e_0e_3-{\lambda^2 \over 2} e_1e_2 & \lambda e_1e_3  &  e_0e_1\cr
 e_3 &e_1e_2+{\lambda^2 \over 2}e_0e_3& 0 & e_0e_1 & -\lambda e_0e_2 \cr
}.$$

\section{powers example} Powers of the maximal ideal of $E$

\def \im {\mathop{\rm im}}

In this section we provide a basic example of the
action of the functors $\LL$ and $\RR$.
Among the most interesting graded $S$-modules are the syzygy modules
that occur in the Koszul complex. We write 
$$
\Omega^i=\coker\bigl[ S\otimes_K\wedge^{i+2}W
\rTo 
S\otimes_K\wedge^{i+1}W\bigr],
$$
where as usual the elements of $W$ have internal degree 1, so that
the generators of $\Omega^i$ have internal degree $i+1$. For example
 $\Omega^{-1}=K$ while $\Omega^0=(W)\subset S$
and $\Omega^{v-1}=S\otimes \wedge^vW$, a free module of rank one 
generated in degree $v$.
The sheafifications of these modules are the exterior
powers of the cotangent bundle on projective space
(see Eisenbud [1995] Section 17.5 for more details.)
In this section we shall show that
under the functors $\LL$ and $\RR$ introduced
in \ref{intro BGG}
the $\Omega^i$
correspond to powers of the maximal ideal $\gm\subset E$. To
make the correspondence completely functorial, we make use of
the $E$-modules $\gm^i\omega_E$, where
$\omega_E=\Hom_K(E,K)$. Recall that $\omega_E$ is a
rank one free $E$-module generated in degree $v$; its
generators may be identified with $\wedge^vW$.

\theorem{powers theorem} 
The minimal $S$-free resolution
of $\Omega^i$ is $\LL(\omega_E/\gm^{v-i}\omega_E)$;
the minimal $E$-injective
resolution of $\omega_E/\gm^{v-i}\omega_E$ is $\RR(\Omega^i)$.

Since $\Omega^i$ is generated in degree $i+1$, 
the complex $\RR(\Omega^i)$
begins in cohomological degree $i+1$, and we regard
$\omega_E/\gm^{v-i}\omega_E$ as concentrated in cohomological 
degree $i+1$.

\proof The complex $\LL(\omega_E)$ is the Koszul complex
over $S$, so $\LL(\omega_E/\gm^{v-i}\omega_E)$ is the truncation
$$
0\rTo S\otimes \wedge^vW \rTo\cdots\rTo S\otimes \wedge^{i+1}W.
$$
which is the resolution of $\Omega^i$, proving the first statement.
The second statement follows from \ref{reciprocity}.\Box

Since the $K$-dual of a minimal $E$-injective resolution is a 
minimal $E$-free resolution, we may immediately derive the
free resolution of
$$
\gm^{i+1}=\Hom_E(\omega_E/\gm^{v-i}\omega_E, \omega_E)
=\Hom_K(\omega_E/\gm^{v-i}\omega_E, K).
$$

\corollary{free res of powers}
The minimal $E$-free resolution of $\gm^{j}$ is
$$
\Hom_K(\RR(\Omega^{j-1}), K).
$$

These resolutions can be made quite explicit using 
the Schur functors $\wedge^i_j$ associated to ``hook'' diagrams
(see for example
Buchsbaum and Eisenbud [1975] or
Akin, Buchsbaum, Weyman [1985]).
We may define
$\wedge^i_j$ (called $L^i_j$ by Buchsbaum and Eisenbud) 
by the formula
$$
\wedge^i_j(W)=\im\bigl[\wedge^iW\otimes_K\sym_{j-1}W\rTo 
                         \wedge^{i-1}W\otimes_K\sym_{j}W\bigr].
$$
Note that 
$$
\wedge^i_j(W)=\cases{0 & if $i<1$ or $j<1$\cr
                     \wedge^iW & if $j=1$\cr 
                     \Sym_jW & if $i=1$
                      }.
$$
Buchsbaum and Eisenbud use these functors to 
give (among other things) 
an $GL(W)$-equivariant resolution 
$$
\cdots\rTo S\otimes_K\wedge^2_j(W)
\rTo S\otimes_K\wedge^1_j(W) 
\rTo (W)^j 
\rTo 0
$$
of the $j^\th$ power
$(W)^j$ of the maximal ideal of $S$.
The $\wedge^i_j$ also provide the terms in the resolutions above:

\corollary{explicit exterior powers} 
The minimal free resolution of $\gm^i$ has the form
$$
\cdots\rTo
E\otimes(\wedge^i_2 W)^*\rTo
E\otimes(\wedge^i_1 W)^*\rTo
\gm^i
\rTo 0
$$
The minimal injective resolution of 
$\omega_E/\gm^{v-i}\omega_E$ has the form
$$
\Hom_K(E, \wedge^{i+1}_{1}(W))\rTo
\Hom_K(E, \wedge^{i+1}_{2}(W))\rTo
\cdots .
$$

\proof From the exactness of the Koszul complex we see that 
$
(\Omega^i)_j = \wedge^{i+1}_{j-i}W, 
$
so the second statement follows from \ref{powers theorem}.
The first statement follows similarly from 
\ref{free res of powers}. \Box

Using the exact sequence
$$
0\to \gm^{v-i}\omega_E\to \omega_E\to \omega_E/\gm^{v-i}\omega_E\to 0
$$
we may paste together the injective and free resolutions considered above
into the Tate resolution $\TT(\Omega^i_\P)$.

\corollary{explicit tate}
There is an exact sequence $\TT(\Omega^i_\P)$
$$\eqalign{
\cdots\rTo
\Hom_K(E,\wedge^{i+1}_1&W)\rTo\Hom_K(E,\wedge^{i+1}_2W)\cr
                        \rTo &\Hom_K(E,K)\rTo           \cr
\Hom_K(E,(\wedge^{v-i}_1&W)^*)\rTo\Hom_K(E,(\wedge^{v-i}_2W)^*)
\rTo\cdots
}$$
where $\Hom_K(E,K)=\omega_E$ is the term in cohomological degree $i$.
\Box

The following well-known result, which we will use
for the proof of \ref{beilinson-thm 2}, now follows from \ref{explicit tate}
by inspection:

\proposition{homology of omegas}
In the range $ 0 \le j \le v-1$  or  $1 \le q \le v-2$    
$$\H^q(\O_\P(-j)\otimes\Omega_\P^p(p)) = 
\cases 
{K, & if p=q=j \cr  
 0, & otherwise. \cr 
}
$$

\proof Writing the ranks of the free modules in the Tate resolution 
for $\Omega_\P^p$ in 
Macaulay notation we find 
$$
\matrix{
& {p+2 \choose 2}{v+2 \choose v-p} & {p+1 \choose 1}{v+1 \choose v-p} & {v \choose v- p} & . & . & . & . \cr
&. & . & . & . & . & . & . \cr
&. & . & . & 1 & . & . & . \cr
&. & . & . & . & . & . & . \cr
&. & . & . & . & . & . & . \cr
&. & . & . & . & {v \choose p-1} & {v-p+2 \choose 1}{v+1 \choose p-1} & {v-p+3 \choose 2} {v+2 \choose p-1}   \cr
}$$
with the rank $1$ module sitting in homological degree $p$ and the out-going
and in-going
 map from it given by bases of the forms in $\wedge^p V$ and
$\wedge^{v-p+1} V$ respectively. (The notation here imitates
that of the Macaulay command ``betti"; it represents the ranks
and degrees of free modules in a
complex with arrows pointing to the left.)
\Box

For the twist $\Omega^p(p)$ we get the rank 1 module in the Tate resolution
in homological degree 0. For their $\Hom$'s we 
prove the following observation of
Beilinson ([1978] Lemma 2) that will play a major role
in \ref{beilinson}.

\proposition{hom of omega example} 
If $\Omega^i(i)$ are the $S$-modules defined in
section \ref{powers example} and $0\leq i,j < v$ then
$$
\hom_S(\Omega^i(i), \Omega^j(j))=
\wedge^{i-j}V = \hom_E(\omega_E(i),\omega_E(j))
$$ 
where in each case $\hom$ denotes the (degree 0) homomorphisms;
for other values of $i,j$ the left hand side is 0.
The product of
homomorphisms corresponds to the product in $\wedge V$.

\proof The modules $\Omega^i(i)$ are 0 for $i<0$ and $i\geq v$.
For $0\leq i<v$  they
have linear resolution, so we may apply \ref{reciprocity}.
As they are 0 in degrees $<1$ and generated in degree 1,
we have 
$\H^1\RR(\Omega^i(i))=\omega_E(i)/\gm^{v-i}\omega_E(i)$ if $v>i$
by \ref{powers theorem}. 
For $0\leq i,j <v$  the maps
$\omega_E(i)/\gm^{v-i}\omega_E(i)\to \omega_E(j)/\gm^{v-j}\omega_E(j)$
are the same as the maps $\omega_E(i)\to\omega_E(j)$. Since $\omega_E$ is 
a rank one free $E$-module, these may be identified with elements of
$E_{j-i}=\wedge^{i-j}V$. \Box

\section{beilinson} Beilinson's Monads

There are two main results. The first says that given 
a sheaf $\F$ on a projective space $\P=\P(W)$ there is a complex
$$
\B:\qquad \dots\rTo B^{-1}\rTo B^0 \rTo B^1\rTo\dots
$$
with
$$
B^e=\oplus_j\H^{j}(\F(e-j))\otimes\Omega_\P^{j-e}(j-e)
$$
such that $\B$ is exact except at $B^0$ and 
the homology at $B^0$ is $\F$.

We show that
the complex $\B$ may be obtained by applying a certain functor to the 
Tate resolution
$\TT(\F)$ over $E$. The second main result gives another monad;
we will deduce it from the first.

Given any graded free complex $\TT$ over $E$ we may
write each module of $\TT$ as a direct sum of copies of 
$\omega_E(i)=\Hom_K(E, K(i))$ with varying $i$. We
define
$\Omega(\TT)$ to be the complex of sheaves on $\P$
obtained by replacing each 
summand $\omega_E(i)$ by the sheaf $\Omega_\P^i(i)$ and using the 
isomorphism of Hom in \ref{hom of omega} to provide the maps.

\theorem{beilinson-thm} If $\F$ is a 
coherent sheaf on $\P(W)$ with associated Tate resolution
$\TT(\F)$,
then the only homology of  $\Omega(\TT(\F)))$ is 
in cohomological degree 0, and is isomorphic to $\F$,

\proof To simplify the notation we set $\TT =\TT (\F)$,
and we let $\overline{\TT}$ be $\TT$
modulo the elements of internal degree
$\geq 0$.  
Let $\L$ be the double complex of sheaves that arises by
sheafifying the double complex of $S$-modules used to 
construct the complex $\LL(\TT)$; that is, if $T^e$ is
the component of $\TT$ of cohomological degree $e$, and
$T^e_j$ is its component of internal degree $j$, then the 
double complex
$\L$ has the form
$$\L:\qquad
\diagram
      &     &\uTo               &     &\uTo\cr
\cdots&\rTo &T^{e}_{j+1}\otimes_K\O(j+1)&\rTo &T^{e+1}_{j+1}\otimes_K\O(j+1)&\rTo&\cdots\cr
      &     &\uTo               &     &\uTo\cr
\cdots&\rTo &T^e_j\otimes_K\O(j)&\rTo &T^{e+1}_{j+1}\otimes_K\O(j)&\rTo&\cdots\cr
      &     &\uTo               &     &\uTo\cr
\enddiagram.
$$
Since $\TT$ is exact, the rows are exact; since the columns are
direct sums of sheafified Koszul complexes over $S$, they are
exact as well.

Choose an integer $f>>0$ (greater than the regularity of 
$\F$ will be sufficient) and let $\L'$ be the double complex
obtained from $\L$ by taking only those terms $T^e_j\otimes_K \O(j)$
with $e<f$ and $j>0$. If $e<<0$ then $T^e$ is generated in negative
degrees, so the double complex $\L'$ is finite, and is exact
except at the right ($e=f-1$) and at j=1. An
easy spectral sequence argument shows that the 
complex obtained as the vertical homology of $\L'$
has the same homology as the complex obtained as the
horizontal homology of $\L'$.

If we write 
$T^e$ as a sum of copies of $\omega_E(i)$ for various $i$,
then the $e^\th$ column of $\L'$ is correspondingly a sum of 
copies of the sheafification of 
$\LL(\omega_E(i)/\gm^{v-i}\omega_E(i)$. As in \ref{powers theorem},
the
vertical homology of this column is correspondingly a sum of 
copies of $\Omega_\P^i(i)$; that is, it is $\Omega(T^e)$. Thus
the complex obtained as the vertical homology of $\L'$ is
$\Omega(\TT)$.

As $e$ goes to infinity the degrees of the generators of $T^e$ 
become more and more positive; thus for $e$ large
the $e^\th$  column of $\L'$ is the same as that of $\L$, that is,
it is $\LL(T^e)$.
Since $f>>0$ the horizontal homology of $\L'$ is 
the sheafification of $\LL(H)$, where $H$ is
the homology of $\TT_{<f}$.
As $\TT$ is exact, $H$ may also be written as the homology of 
$\TT_{\geq f}$. Taking $f>\reg\F$ and using \ref{reciprocity},
we see that $\LL(H)$ is a free resolution of the module
$\oplus_{e\geq f}\H^0(\F(e))$, whose sheafification is $\F$,
as required.\Box

\corollary{beilinson cor 2} The map in 
the complex $\Omega(\TT (\F))$ corresponding to 
$$
\H^j(\F(j-i))\otimes\Omega_\P^{i-j}(i-j)\rTo
\H^j(\F(j-i+1))\otimes\Omega_\P^{i-j-1}(i-j-1)
$$
corresponds to the multiplication map 
$W\otimes_K \H^j(\F(j-i))\rTo \H^j(\F(j-i+1))$.

\proof This follows from \ref{sheaf cohomology}, since we have identified
not only the modules but the maps in the linear strands of the 
resolution.
\Box

\corollary{beilinson cor 1} The maps in the complex
$\Omega(\TT (\F))$ correspond to the maps in the complex
$\Omega(\TT (\F(1)))$ under the natural correspondence
$$
\hom_\P(\Omega_\P^i(i), \Omega_\P^j(j))=\wedge^{i-j}V
=\hom_\P(\Omega_\P^{i+1}(i+1), \Omega_\P^{j+1}(j+1))
$$
whenever $0\leq i, i+1, j, j+1 <v$.

\proof The Tate resolution $\TT (\F(1))$ is obtained by simply
shifting $\TT (\F)$.\Box

Beilinson's second theorem gives similarly a description of $\F$
as the homology of a complex of sheaves whose terms are direct sums
of  line bundles $\O_\P,\ldots,\O_\P(-n)$ with $n=v-1$. 
Rather than deduce it with a proof
parallel to that of \ref{beilinson-thm} 
we will deduce it from \ref{beilinson-thm}.

\theorem{beilinson-thm 2} Set
$n=\dim\ \P(W)$. For every coherent sheaf $\F$ on $\P(W)$
there is a minimal finite complex
$$\B: \ldots \rTo B^{-1} \rTo B^0 \rTo B^1 \rTo
\ldots$$
with
$$
B^j = \oplus_{q=0}^n \H^{j+q}(\F \otimes \Omega_\P^{q}(q))
\otimes \O_\P(-q),
$$ 
whose only homology is $\H^0(\B) \cong \F$.

\proof 
Consider the Beilinson monad $\Omega(\F(-1)) \tensor \O_\P(1)$.
Each term in the complex is a direct sum of sheaves $\Omega_\P^i(i+1)$.
We consider the double complex $\C=\C(\F)$ obtained by replacing each 
summand $\Omega_\P^i(i+1)$ by its resolution 
$$\eqalign{
\LL(\omega_E&(i+1)/\gm^{v-i}\omega_E(i+1)):\cr
&
\quad
0
\rTo 
\wedge^vW\otimes_K \O_\P(i-n)
\rTo 
\dots
\rTo
\wedge^{i+1}W\otimes_K \O_\P
\rTo 
0}
$$
and by replacing each component of a differential in $\Omega(\F(-1))
\otimes \O_\P(1)$
  by the 
morphism of complexes 
$$
\LL(\omega_E(i+1)/\gm^{v-i}\omega_E(i+1))\to 
\LL(\omega_E(j+1)/\gm^{v-i}\omega_E(j+1))
$$
defined by the corresponding map $\omega_E(i)\to \omega_E(j)$ in $\TT(\F)$.
The terms of the double complex $\C(\F)$ are sums of line bundles
$\O_\P,\ldots,\O_\P(-n)$. The maps in the rows
are of degree 0, while the columns are truncated
Koszul complexes. The vertical homology
of $\C(\F)$ is the Beilinson monad 
$\Omega(\F(-1)) \tensor \O_\P(1)$, which has as its only homology $\F$ in
degree 0. So the total complex of $\C(\F)$ has the same homology

We now make $\C(\F)$ {\it minimal,\/} by successively factoring out
subcomplexes which are direct summands and have the form
$$
\cdots 0\rTo \O_\P(i)\rTo^a \O_\P(i)\rTo 0 \rTo\cdots
$$
where $a$ is nonzero, and is thus an isomorphism. We eventually
arrive at a complex
$\B(\F)$ which is minimal in the sense that it contains no such
summands. This process does not change the homology, so the homology
of $\B$ is still $\F$, concentrated in cohomological degree 0. The
proof of \ref{beilinson-thm 2} is completed by the following:

\proposition{beilinson-uniqueness} Let $\B$ be a complex of 
sheaves on $\P=\P(W)$ whose terms are sums of line bundles
$\O_\P,\ldots,\O_\P(-n)$, where $n=v-1=\dim\ \P$, and which is minimal
in the sense above. If the homology of $\B$ is 
the sheaf $\F$, concentrated in cohomological degree 0, then
the $j^\th$ term $B^j$ of $\B$ is 
$$
B^j = \oplus_{q=0}^n \H^{j+q}(\F \otimes \Omega_\P^{q}(q))
\otimes \O_\P(-q).
$$

\proof
By \ref{homology of omegas} we can ``detect'' copies of
$\O_\P(-q)$ in $B^j$ by tensoring with $\Omega_P^q(q)$
and taking cohomology. A diagram chase completes the proof.
\Box

\section{Examples} Examples

\example{rat normal curve} Let $C \subset \P^d$ be the rational normal
curve parametrized by $(s:t) \mapsto (s^d:s^{d-1}t:\ldots:t^d)$. We consider
 the line bundles $\L_k$ on $C$ associated to 
$\oplus_{m=0}^\infty H^0(\P^1,\O(k+md))$ for $k=-1,\ldots,d-2$. The Tate 
resolution
$\TT(\L_k)$ has betti numbers
$$
\matrix{ * \; * & 2d+k+1 & d+k+1 & k+1   & .     & .      &\ldots \cr
	 \ldots & .      & .      & d-k-1 & 2d-k-1 & 3d-k-1  &*\;* \cr
}
$$
The $(k+1+d-k-1) \times (2d-k-1)$ middle matrix and the matrices surrounding it
 have in case $d=4$ and $k=1$ the following shapes:
$$\scriptscriptstyle{\pmatrix{
0&e_0&e_0e_2&e_0e_1 \cr
e_0&e_1&e_1e_2+e_0e_3&e_0e_2 \cr
e_1&e_2&e_1e_3+e_0e_4&e_0e_3 \cr
e_2&e_3&e_1e_4&e_0e_4\cr
e_3&e_4&0&0\cr
e_4&0&0&0\cr},
\pmatrix{ 
0&0& e_0e_4& e_1e_4        &   e_2e_4 & e_3e_4 \cr
0&0& e_0e_3& e_1e_3+e_0e_4 & e_2e_3+e_1e_4 & e_2e_4 \cr
0&e_0&e_1&e_2&e_3&e_4\cr
e_0&e_1&e_2&e_3&e_4&0\cr	
 }}
$$
and
$$\pmatrix{
e_0&e_1&e_2&e_3&e_4&0&0&0&0&0\cr
0&e_0&e_1&e_2&e_3&e_4&0&0&0&0\cr
0&0&e_0&e_1&e_2&e_3&e_4&0&0&0\cr
0&0&0&e_0&e_1&e_2&e_3&e_4&0&0\cr
0&0&0&0&e_0&e_1&e_2&e_3&e_4&0\cr
0&0&0&0&0&e_0&e_1&e_2&e_3&e_4\cr
}$$
All other matrices look similar to the last one.

Notices that in case $k=-1$ we obtain a $(d-1)\times (d-1)$ symmetric 
matrix of 2-forms:
$$\pmatrix{
e_0e_1&e_0e_2&e_0e_3&e_0e_4\cr
e_0e_2&e_1e_2+e_0e_3&e_1e_3+e_0e_4& e_1e_4\cr
e_0e_3&e_1e_3+e_0e_4&e_2e_3+e_1e_4&e_2e_4 \cr
e_0e_4&e_1e_4&e_2e_4 & e_3e_4\cr
}.
$$
If we interpret  2-forms as coordinate functions 
$$e_{ij}=e_ie_j=e_i\wedge e_j \in \H^0(\GG(W,2),\O(1)) \cong
 \H^0(\P(\Lambda^2 V),\O(1))$$
 on the Grassmanian
of codimension 2 linear subspaces in $\P(W)$, then the  determinant of the 
matrix above defines the Chow point of $C\subset \P^d$, which is by definition the 
hypersurface
 $\{\P^{d-2} \in \GG(W,2) | \P^{d-2} \cap C \ne \emptyset \}$.
We will prove this in a general setting in a forthcoming paper.

\example{Horrocks-Mumford} A famous Beilinson monad is due to Horrocks and 
Mumford:
Consider for $\P^4$ the Tate resolution $\TT(\varphi)$ of the matrix
$$
\varphi = \pmatrix{ 
e_1e_4 & e_2e_0 &e_3e_1 & e_4e_2 & e_0e_3 \cr
e_2e_3 & e_3e_4 &e_4e_0 & e_0e_1 & e_1e_2 \cr 
}.$$
By straight computation we find the betti numbers 
$$\matrix{ 
*\; *&100&35 & 4 & . & . & . & . & . & . & . & . & \ldots \cr
\ldots & . & 2 & 10&10 & 5 & . & . & . & . & . & . &\ldots \cr
\ldots & . &.  & . & . & . & 2 & . & . & . & . & . &\ldots \cr
\ldots & . &.  & . & . & . & . & 5 & 10&10 & 2 & . &\ldots \cr
\ldots & . &.  & . & . & . & . & . & . & 4 & 35&100&*\;* \cr
}.$$

To deduce that this Tate resolution comes from a sheaf we use:

\lemma{row bound} Let $\TT$ be a Tate resolution over $E$. Suppose 
 that  $(T^0 \tensor K)_j = 0$  for all $j < r$. 
Then 
$(T^l \tensor K)_m = 0$  if 
$ l>0$  and  $ m-l <r$, or if
$l<0$ and $m-l<r-v+1$.

\proof  The first vanishing follows, because  $\Hom(\TT,E)$ is also a minimal
complex. For the second we note for $P=\ker(T^0 \to T^1)$
that $P_j=0$ holds for all $j<r-v$ by our assumption. By \ref{explicit exterior powers} 
$(T^l\tensor K)=\Tor^E_{-l-1}(P,K)=(P\tensor_K (\Sym_{-l-1} W)^*)$. So this group vanishes in all degree $m<r-v+l+1$.
\Box

\noindent{\bf\ref{Horrocks-Mumford} Continued:} 
By applying  \ref{row bound} to a shifts of  $\TT(\varphi)$ and 
$\Hom(\TT(\varphi),E)$ we
see that the $\TT(\varphi)$ has  terms only in the indicated  
range of rows. So $\TT(\varphi)$ is the Tate resolution of some sheaf $\F$.
Moreover $\F$ is a bundle, since the middle cohomology has only finitely many
terms. The $4^{\th}$ difference function of  $\chi(\F(m))$
has constant value $2$. So
$\F$ has rank $2$. It is the famous  bundle on $\P^4$
discovered by Horrocks and Mumford [1973]. In Decker and Schreyer [1986] it
is proved that any stable rank $2$ vector bundle on $\P^4$ with the same 
Chern classes equals $\F$ up to projectivities.

\references
\frenchspacing
\parindent=0pt

\noindent
K.~Akin, D.~A.~Buchsbaum, and J.~Weyman: Schur functors and Schur
complexes. Adv. in Math.~44 (1982) 207--278.

\medskip
\noindent
V.~Ancona and G.~Ottaviani: An introduction to the derived
categories and the theorem of Beilinson. Atti Accademia 
Peloritana dei Pericolanti, Classe I de Scienze Fis. Mat. et Nat.
LXVII (1989) 99--110.

\medskip
\noindent
A.~Aramova, L.A.~Avramov, J.~Herzog: 
Resolutions of monomial ideals and cohomology over exterior algebras.
Trans. Amer. Math. Soc. 352 (2000)579--594.

\medskip
\noindent
A. Beilinson: Coherent sheaves on $\P^n$ and problems
of linear algebra. 
Funct. Anal. and its Appl. 12 (1978) 214--216.
(Trans. from Funkz. Anal. i. Ego Priloz 12 (1978) 68--69.)

\medskip
\noindent
Bernstein, Gel'fand, Gel'fand:
Algebraic bundles on $\P^n$ and problems of linear algebra.
Funct. Anal. and its Appl. 12 (1978) 212--214.
(Trans. from Funkz. Anal. i. Ego Priloz 12 (1978) 66--67.)

\medskip
\noindent
D.~A.~Buchsbaum and D.~Eisenbud: Generic Free resolutions and
a family of generically perfect ideals. Adv.~Math.~18 (1975) 245--301.
\medskip
\noindent

R.-O. Buchweitz: Appendix to 
Cohen-Macaulay modules on quadrics, by R.-O. Buchweitz, D. Eisenbud,
and J. Herzog. In {\sl Singularities,
representation of algebras, and vector bundles (Lambrecht, 1985),\/}
Springer-Verlag Lecture Notes in Math. 1273 (1987) 96--116.

\medskip
\noindent
W.~Decker and F.-O.~Schreyer: On the uniqueness of the 
Horrocks-Mumford bundle. Math. Ann. 273 (1986) 415--443.

\medskip
\noindent
W.~Decker: Stable rank 2 bundles with Chern-classes
$c_1=-1,\ c_2=4$. Math. Ann. 275 (1986) 481--500.

\medskip
\noindent
D. Eisenbud: {\sl Commutative Algebra with a View Toward 
Algebraic Geometry.} Springer Verlag, 1995.

\medskip
\noindent
D. Eisenbud and S. Goto: Linear free resolutions and minimal
multiplicity.  J. Alg. 88 (1984) 89--133.

\medskip
\noindent
D. Eisenbud and S. Popescu: 
Gale Duality and Free Resolutions of Ideals of Points.
Invent. Math. 136 (1999) 419--449.

\medskip
\noindent
D. Eisenbud, S. Popescu, and S. Yuzvinsky: Hyperplane
arrangements and resolutions of monomial ideals over an
exterior algebra. Preprint (1999).

\medskip
\noindent
S. I. Gel'fand and Yu. I. Manin: {\it Methods of Homological Algebra.\/}
Springer-Verlag, New York 1996.

\medskip
\noindent
S. I. Gel'fand: Sheaves on $\P^n$ and problems in linear algebra.
Appendix to the Russian edition
of K.~Okonek, M.~Schneider, and K.~Spindler. Mir, Moscow, 1984.

\medskip
\noindent
M. Green: The Eisenbud-Koh-Stillman Conjecture on Linear Syzygies.
Invent. Math. 136 (1999) 411--418.

\medskip
\noindent
D.~Grayson and M.~Stillman: Macaulay2.\hfill\break
{\tt http://www.math.uiuc.edu/Macaulay2/}.

\medskip
\noindent
A. Grothendieck and J. Dieudonn\'e: 
{\it \'El\'ements de la G\'eometrie Alg\'ebrique IV:
\'Etude locale des sch\'emas et de morphismes de sch\'emas}
(deuxi\`eme partie). Publ. Math. de l'I.H.E.S. 24 (1965).

\medskip
\noindent
D.~Happel: {\it Triangulated Categories In The Representation Theory
Of Finite-Dimensional Algebras.\/} Lect.~Notes of the London
Math.~Soc.~119, 1988.

\medskip
\noindent
R. Hartshorne: {\it Algebraic geometry \/}, Springer Verlag, 1977.

\medskip
\noindent
J. Herzog, T. R\"omer: Resolutions of modules over the exterior algebra,
working notes, 1999.

\medskip
\noindent
G. Horrocks, D. Mumford: A rank $2$ vector bundle on $\P^4$ with $15,000$
symmetries. Topology 12 (1973) 63--81.

\medskip
\noindent
 M. M. Kapranov:
On the derived categories of coherent sheaves on some homogeneous spaces. 
Invent. Math. 92 (1988) 479--508.

\medskip
\noindent
 M. M. Kapranov:
On the derived category and $K$-functor 
of coherent sheaves on intersections of quadrics. (Russian) Izv. Akad.
Nauk SSSR Ser. Mat. 52 (1988) 186--199; translation in Math.
USSR-Izv. 32 (1989)  191--204.

\medskip
\noindent
Ch. Okonek, M. Schneider, and H. Spindler: {\sl Vector Bundles on
Complex Projective Spaces.\/} Birkh\"auser, Boston 1980.
\medskip
\noindent

 D. O. Orlov:
Projective bundles, monoidal transformations, 
and derived categories of coherent sheaves. (Russian) Izv. Ross. Akad. Nauk
Ser. Mat. 56 (1992) 852--862; translation in Russian Acad. Sci. Izv. Math.
41 (1993) 133--141.

\medskip
\noindent
S.B. Priddy: Koszul resolutions. Trans. Amer. Math. Soc. 152
(1970) 39--60. 

\medskip
\noindent
R.~G.~Swan: $K$-theory of quadric hypersurfaces. Ann. of Math.
122 (1985)113--153.

\bigskip
\vbox{\noindent Author Addresses:
\smallskip
\noindent{David Eisenbud}\par
\noindent{Department of Mathematics, University of California, Berkeley,
Berkeley CA 94720}\par
\noindent{eisenbud@math.berkeley.edu}
\smallskip
\noindent{Frank-Olaf Schreyer}\par
\noindent{FB Mathematik, Universit\"at Bayreuth
D-95440 Bayreuth, Germany\par}
\noindent{schreyer@btm8x5.mat.uni-bayreuth.de}\par
}

\end